\documentclass[journal]{IEEEtran}
\usepackage{cite}
\usepackage{amsmath}
\usepackage{amssymb}
\newtheorem{theorem}{Theorem}
\newtheorem{proposition}{Proposition}
\newtheorem{lemma}{Lemma}
\newcommand*{\QEDB}{\hfill\ensuremath{\square}}%
\usepackage{footnote}
\makesavenoteenv{tabular}
\usepackage{algorithm}
\usepackage{algorithmic}
\newtheorem{remark}{Remark}
\usepackage{mleftright}
\usepackage{epstopdf}

\usepackage{tikz}
\usetikzlibrary{shapes,arrows,fit,positioning,shadows,calc}
\usetikzlibrary{plotmarks}
\usetikzlibrary{decorations.pathreplacing}
\usetikzlibrary{patterns}
\usetikzlibrary{automata}
\usepackage{pgfplots}
\pgfplotsset{compat=newest}
\usepackage{adjustbox}
\usepackage{caption}
\begin{document}
%
\title{Subspace-Orbit Randomized Decomposition for Low-rank Matrix Approximation}

\author{\IEEEauthorblockN{Maboud F. Kaloorazi$^{\dagger}$, and 
Rodrigo C. de Lamare$^{\dagger,\ddagger}$} \\
\vspace{5mm}
\IEEEauthorblockA{$^\dagger$ Centre for Telecommunications Studies (CETUC)}\\
Pontifical Catholic University of Rio de Janeiro, Brazil \\
$^\ddagger$ Department of Electronics, University of York, United Kingdom \\
E-mail: \texttt{\{kaloorazi,delamare\}@cetuc.puc-rio.br}}

%
%

\IEEEtitleabstractindextext{%
\begin{abstract}
An efficient, accurate and reliable approximation of a matrix by one of lower 
rank is a fundamental task in numerical linear algebra and signal processing applications. In this paper, we introduce a new matrix decomposition approach 
termed Subspace-Orbit Randomized singular value decomposition (SOR-SVD), which 
makes use of random sampling techniques to give an approximation to a low-rank matrix. Given a large and dense data matrix of size $m\times n$ with numerical rank $k$, where $k \ll \text{min} \{m,n\}$, the algorithm requires a few passes through data, and can be computed in $O(mnk)$ floating-point operations. Moreover, the SOR-SVD algorithm can utilize advanced computer architectures, and, as a result, it can be optimized for maximum efficiency. The SOR-SVD algorithm is simple, accurate, and provably correct, and outperforms previously reported techniques in terms of accuracy and efficiency. Our numerical experiments support these claims.
\end{abstract}

\begin{IEEEkeywords}
Matrix computations, dimension reduction, low-rank approximation, randomized algorithms.
\end{IEEEkeywords}}

\maketitle

\IEEEdisplaynontitleabstractindextext

%
\IEEEpeerreviewmaketitle

\section{Introduction}\label{sec:intr}

\IEEEPARstart{C}{omputing} a low-rank approximation of a given data matrix, i.e., approximating the matrix by one of lower rank, is a fundamental task in signal processing and its applications. Such a compact representation of a matrix can provide a significant reduction in memory requirements, and more importantly, computational costs when the computational cost scales according to a high-degree polynomial with the dimensionality. Matrices with low-rank structures arise frequently in numerous applications such as latent variable graphical modeling, \cite{Chandrasekaran12}, ranking and collaborative filtering, \cite{Srebro2005}, background subtraction \cite{WPMGR2009, BZ2014, KaDeDSP17, RahmaniAtiaHighP17}, system identification \cite{FazelPST13}, IP network anomaly detection \cite{MMG2013,KaDeICASSP17}, 
subspace clustering \cite{SoltanolkotabiEC2014, RahmaniAtiaCoP17, Oh2017, ShenMG2017}, biometrics \cite{VictorBS02, WrightFace09}, sensor and multichannel signal 
processing \cite{DeSa2009}, statistical process control and multidimensional fault 
identification \cite{Jackson91, Dunia1998}, DNA microarray data \cite{Troyanskaya01missingvalue}, and quantum state tomography \cite{Gross2010}. 

\stepcounter{footnote}\footnotetext {A primary version of this paper has been submitted to the 26th European Signal Processing Conference (EUSIPCO 2018).}

Traditional algorithms for computing a low-rank approximation of a matrix, such as singular value decomposition (SVD) \cite{GolubVanLoan96} and the rank-revealing QR decomposition \cite{Chan87, GuEisenstat96} are computationally prohibitive for large data sets. Moreover, standard techniques for their computation are challenging to parallelize in order to take advantage of modern computer architectures \cite{Demmel97, MarQH17}. Recently developed low-rank approximation algorithms based random sampling techniques, however, have been shown to be surprisingly efficient, accurate and robust, and are known to outperform the traditional algorithms in many practical situations \cite{FriezeKVS04, DrineasKM06, Sarlos06, Rokhlin09, HMT2009,Gu2015}. The power of randomized methods lies in the fact they can be optimized for maximum efficiency on modern architectures.
\begin{itemize}
\item [] \textit{Our Contributions}
\end{itemize}

Motivated by recent developments, this paper proposes a new randomized decompositional approach called subspace-orbit randomized singular value decomposition (SOR-SVD) to compute the low-rank approximation of a matrix. 
Given a large and dense rank-$k$ matrix of size $m\times n$, the SOR-SVD 
requires a few passes over the data, and can be computed in $O(mnk)$ 
floating-point operations (flops). The main operations of the algorithm 
involve matrix-matrix multiplication and the QR decomposition. This, due 
to recently developed Communication-Avoiding QR (CAQR) algorithms 
\cite{DemmGHL12}, which are optimal in terms of communication costs,  i.e., 
data movement either between different levels of a memory hierarchy or 
between processors, allows the algorithm to be optimized for peak machine performance on modern computational platforms.

We provide theoretical error bounds, i.e., lower bounds on singular values 
and upper bounds on the low-rank approximation, for the SOR-SVD algorithm, 
and experimentally show that the low-rank error bounds provided are 
empirically sharp for one class of matrices considered.
   
We also apply the SOR-SVD to solve the RPCA problem \cite{WPMGR2009,CSPW2009, 
CLMW2009}, i.e., to decompose a data matrix into a low-rank plus a sparse 
matrix, and study in computer vision applications of background/foreground separation in surveillance video, and shadow and specularity removal from 
face images.

\begin{itemize}
\item [] \textit{Notation}
\end{itemize}

Bold-face upper-case letters are used to denote matrices. Given a matrix $\bf A$, ${\|{\bf A}\|_1}$, ${\|{\bf A}\|_2}$, ${\|{\bf A}\|_F}$, ${\|{\bf A}\|_*}$ denote the $\ell_1$-norm, the spectral norm, the Frobenius norm, the nuclear norm, respectively. $\sigma_j(\bf A)$ denotes the $j$-th largest singular value of 
$\bf A$, and the numerical range of $\bf A$ is denoted by $\mathcal{R}({\bf A})$. 
The symbol $\mathbb{E}$ denotes expected value with respect to random variables. 
Given a random variable $\Omega$, $\mathbb{E}_\Omega$ denotes expectation with respect to the randomness in $\Omega$, and the dagger $\dagger$ denotes the 
Moore-Penrose pseudo-inverse.

We structure the remainder of this paper as follows. In Section \ref{secRelatW}, we introduce the mathematical model of the data and discuss prior works on randomized algorithms. In Section \ref{secSOR}, we describe our proposed approach, which also includes a variant that uses the power iteration scheme, in detail. Section \ref{secAnalysis} presents our theoretical analysis. In Section \ref{secNumExp}, we present and discuss our numerical experimental results, and our conclusion remarks are given in Section \ref{secCon}.

\section{Mathematical Model and Prior Works}
\label{secRelatW}
Given an input matrix ${\bf A} \in \mathbb R^{m \times n}$, where $m \ge n$, 
with numerical rank $k$, its singular value decomposition (SVD) \cite{GolubVanLoan96} is defined as follows: 
\begin{equation}
\begin{aligned}
{\bf A} = & {\bf U}_{m\times n}{\bf \Sigma}_{n\times n}{\bf V}_{n\times n}^T \\ = & \begin{bmatrix} {{\bf U}_k \quad {\bf U}_0} \end{bmatrix}
  \begin{bmatrix}
       {\bf \Sigma}_k & 0  \\
       0 & {\bf \Sigma}_0
  \end{bmatrix}\begin{bmatrix}{{\bf V}_k \quad {\bf V}_0} \end{bmatrix}^T,
\label{eqSVD}
\end{aligned}
\end{equation}
where ${\bf U}_k \in \mathbb R^{m \times k}$, ${\bf U}_0 \in \mathbb R^{m \times n-k}$ have orthonormal columns, spanning the range of $\bf A$ and the null space of ${\bf A}^T$, respectively, ${\bf \Sigma}_k \in \mathbb R^{k \times k}$ and 
${\bf \Sigma}_0 \in \mathbb R^{n-k \times n-k}$ are diagonal containing the singular values, i.e., ${\bf \Sigma}_k=\text{diag}(\sigma_1, ..., \sigma_k)$   and ${\bf \Sigma}_0 =\text{diag}(\sigma_{k+1}, ..., \sigma_n)$, and ${\bf V}_k \in \mathbb R^{n \times k}$ and ${\bf V}_0 \in \mathbb R^{n \times n-k}$ have orthonormal columns, spanning the range of ${\bf A}^T$ and the null space of ${\bf A}$, respectively. $\bf A$ can be written as ${\bf A} = {\bf A}_k+{\bf A}_0$, where ${\bf A}_k = {\bf U}_k{\bf \Sigma}_k{\bf V}_k^T$, and ${\bf A}_0 = {\bf U}_0{\bf \Sigma}_0{\bf V}_0^T$. The SVD  constructs the optimal rank-$k$ approximation ${\bf A}_k$ to ${\bf A}$, as stated in the following theorem.

\begin{theorem}
(Eckart and Young \cite{EckartYoung36}, and Mirsky \cite{Mirsky60})
\begin{equation}
\begin{aligned}
& \underset{\text{rank}({\bf B})\le k}{\text{minimize}}
&& \|{\bf A} - {\bf B}\|_2 = \|{\bf A} - {\bf A}_k\|_2 = \sigma_{k+1}.
\end{aligned}
\end{equation}
\begin{equation}
\begin{aligned}
&\underset{\text{rank}({\bf B})\le k}{\text{minimize}}
&&\|{\bf A} - {\bf B}\|_F = \|{\bf A} - {\bf A}_k\|_F = \sqrt{\sum_{j=k+1}^{n}{\sigma_j^2}}.
\end{aligned}
\end{equation}
\label{ThrEckYou}
\end{theorem}

Throughout this paper we focus on the matrix $\bf A$ defined, and $m \ge n> \text{max}\{k, 2\}$. Our results hold for $m < n$, though.

Although, the SVD is numerically stable, accurate and provides detailed 
information on singular subspaces and singular values, it is expensive to 
compute for large data sets. Furthermore, standard techniques for its 
computation are challenging to parallelize in order to utilize modern 
computational platforms \cite{Demmel97, MarQH17}. 
On the other hand, partial SVD based on Krylov subspace methods, such as the 
Lanczos and Arnoldi algorithms, can construct an approximate SVD of an input 
matrix, for instance $\bf A$, best, at a cost $O(mnk)$. However, the partial SVD 
suffers from two drawbacks. First, inherently, it is numerically unstable \cite{CalvettiRS94,GolubVanLoan96, Demmel97}. Second, its methods are 
challenging to parallelize \cite{HMT2009, Gu2015}, which makes it 
unsuitable to apply on modern architectures. 

Recently developed algorithms based on randomization for low-rank approximations  \cite{FriezeKVS04, DrineasKM06, Sarlos06, Rokhlin09, HMT2009, Gu2015,TrYUC16} have attracted considerable attention because \textit{i)} they are computationally efficient, and \textit{ii)} they readily lend themselves to a parallel implementation to exploit modern computational architectures.

Beginning with \cite{FriezeKVS04}, many algorithms have been proposed for 
low-rank matrix approximation. The algorithms in \cite{DrineasKM06, 
DeshpandeV2006, RudelsonV07}, bulit on Frieze et al.'s idea \cite{FriezeKVS04}, first sample columns of an input matrix with a probability proportional to 
either their magnitudes or leverage scores, representing the matrix in a 
compressed form. The submatrix is then used for further computation 
(post-processing step) using deterministic algorithms such as the SVD and 
pivoted QR decomposition \cite{GolubVanLoan96} to obtain the final low-rank approximation. Sarl{\'{o}}s \cite{Sarlos06} proposes a different method based on results of the well-known Johnson-Lindenstrauss (JL) lemma \cite{JL84}. He showed that random linear combinations of rows, i.e., projecting the data matrix onto a structured random subspace, can render a good approximation to a low-rank matrix. 
The works in \cite{NelsonNg2013,ClarWood2017} have further advanced Sarl{\'{o}}s's idea and construct a low-rank approximation based on subspace embedding. Rokhlin et al. \cite{Rokhlin09} proposes to apply a random Gaussian embedding matrix in order to reduce the dimension of the data matrix. The low-rank approximation is then given through computations using the classical techniques on the reduced-sized matrix. In a seminal paper, Halko et al. \cite{HMT2009} propose two algorithms based on 
randomization. (The paper also develops modified versions of each for 
gaining either performance or pass-efficiency). The first algorithm, 
\textit{randomized SVD}, for which the authors provide theoretical analysis 
and extensive numerical experiments, for the matrix ${\bf A}$, and integers 
$k\le \ell<n$ and $q$, is described in Alg. \ref{Alg1}.\looseness-1
\begin{algorithm}
\caption{Randomized SVD (R-SVD)}
\renewcommand{\algorithmicrequire}{\textbf{Input:}}
\begin{algorithmic}[1]
\REQUIRE ~~ 
 Matrix $\ {\bf A} \in \mathbb R^{m \times n}$,
integers $\ell$ and $q$,
\renewcommand{\algorithmicrequire}{\textbf{Output:}}
\REQUIRE ~~ A rank-$\ell$ approximation.
  \STATE Draw a random matrix ${\bf \Omega} \in \mathbb R^{n \times \ell}$;
  \STATE Compute ${\bf Y} = ({\bf A A}^T)^q{\bf A}{\bf \Omega}$;
  \STATE Compute a QR decomposition ${\bf Y} = {\bf Q}{\bf R}$; 
  \STATE Compute ${\bf B} = {\bf Q}^T\bf A$;
  \STATE Compute an SVD ${\bf B} = \widetilde{\bf U} {\bf \Sigma}{\bf V}^T$;
  \STATE ${\bf A} \approx {\bf Q} \widetilde{\bf U}{\bf \Sigma}{\bf V}^T$.
\end{algorithmic}\label{Alg1}
\end{algorithm}

The R-SVD approximates $\bf A$ as follows: \textit{i)} a 
compressed matrix $\bf Y$, through random linear combinations of columns of 
$\bf A$ is formed, \textit{ii)} a QR decomposition is performed on $\bf Y$, 
where the $\bf Q$ factor constructs an approximate basis for $\mathcal{R}
({\bf A})$, \textit{iii)} $\bf A$ is projected onto a subspace spanned by 
columns of $\bf Q$, forming $\bf B$, \textit{iv)} a full SVD of $\bf B$ is 
computed. In Alg. \ref{Alg1}, $q$ is the number of steps of a power method \cite{Rokhlin09, HMT2009}. 

Gu \cite{Gu2015} applies a slightly modified version of the R-SVD algorithm to improve subspace iteration methods, and presents a new error analysis. The second method proposed in \cite[Section 5.5]{HMT2009} is a \textit{single-pass} algorithm, i.e., it requires only one pass through data, to compute a low-rank approximation. For the matrix ${\bf A}$, the decomposition, which we call two-sided randomized SVD (TSR-SVD), is computed as described in Alg. \ref{Alg2}.\looseness-1
\begin{algorithm}
\caption{Two-Sided Randomized SVD (TSR-SVD)}
\renewcommand{\algorithmicrequire}{\textbf{Input:}}
\begin{algorithmic}[1]
\REQUIRE ~~ 
 Matrix $\ {\bf A} \in \mathbb R^{m \times n}$,
an integer $\ell$,
\renewcommand{\algorithmicrequire}{\textbf{Output:}}
\REQUIRE ~~ A rank-$\ell$ approximation.
  \STATE Draw random matrices ${\bf \Psi}_1 \in \mathbb R^{n \times \ell}$ and 
  ${\bf \Psi}_2 \in \mathbb R^{m \times \ell}$;
  \STATE Compute ${\bf Y}_1 = {\bf A}{\bf \Psi}_1$ and ${\bf Y}_2 = {\bf A}^T{\bf \Psi}_2$ in a single pass through $\bf A$;
  \STATE Compute QR decompositions ${\bf Y}_1 = {\bf Q}_1{\bf R}_1$, ${\bf Y}_2 = {\bf Q}_2{\bf R}_2$; 
  \STATE Compute ${\bf B}_\text{approx} = {\bf Q}_1^T{\bf Y}_1({\bf Q}_2^T
  {\bf \Psi}_1)^\dagger $;
  \STATE Compute an SVD ${\bf B}_\text{approx} = \widetilde{\bf U} 
  \widetilde{\bf \Sigma} \widetilde{\bf V}$;
  \STATE ${\bf A} \approx ({\bf Q}_1 \widetilde{\bf U})\widetilde{\bf \Sigma}
  ({\bf Q}_2 \widetilde{\bf V})^T$.
\end{algorithmic}\label{Alg2}
\end{algorithm}

Here ${\bf B}_\text{approx}$ is an approximation to ${\bf B} = {\bf Q}_1^T{\bf A}{\bf Q}_2$. 

The TSR-SVD, however, substantially degrades the quality of the 
approximation, compared to the R-SVD, even for matrices with rapidly 
decaying singular values. The reason behind is, mainly, poor approximate 
basis drawn from the row space of $\bf A$, i.e., ${\bf Q}_2$. Furthermore, 
for a general input matrix, the authors do not provide theoretical bounds, 
neither upper bounds on the low-rank approximation nor lower bounds on the 
singular values, for the TSR-SVD algorithm. This work addresses these issues.

The work in \cite{KaDeDSP17} proposes a rank-revealing decomposition method 
based on randomized sampling techniques; the matrix ${\bf A}$ is compressed, multiplied by orthonormal bases for $\mathcal{R}({\bf A})$ and $\mathcal{R}({\bf A}^T)$ from both sides, columns of the reduced matrix, and accordingly, 
the bases are permuted, the low-rank approximation is given by projecting it 
back to the original space. 

Independent of the TSR-SVD algorithm \cite{HMT2009}, this work 
was developed by drawing inspiration from the randomized rank-revealing 
decomposition algorithm proposed in \cite{KaDeDSP17}. Unlike the R-SVD, 
Alg. \ref{Alg1}, \cite{HMT2009, Gu2015}, which is based on \textit{one-sided} projection, the SOR-SVD algorithm, similar to the TSR-SVD algorithm, 
is based on \textit{two-sided} projection; the data matrix is projected 
onto a lower-dimensional space using orthonormal bases from both sides, 
a truncated SVD of the compressed matrix is computed, the approximate SVD 
is obtained by projecting the data back to the original space. 

We also present a variant of the SOR-SVD algorithm, which employs the power 
method scheme \cite{Rokhlin09, HMT2009} to improve the accuracy of the approximations.

Our analysis was inspired by the excellent work of Gu \cite{Gu2015}. However, 
unlike the one-sided randomized SVD algorithm presented in \cite{Gu2015}, 
we perform an analysis of a two-sided randomized algorithm, i.e., the SOR-SVD, 
and provide theoretical lower bounds on singular values and upper bounds 
on the low-rank approximation. 

\section{Subspace-Orbit Randomized Singular Value Decomposition}
\label{secSOR}
In this section, we present a randomized algorithm termed subspace-orbit 
randomized SVD (SOR-SVD), which computes the \textit{fixed-rank} approximation 
of a data matrix. We also present a version of the SOR-SVD with power 
iteration, which improves the performance of the algorithm at an extra 
computational cost.

Given the matrix ${\bf A}$, and an integer $k\le \ell< n$, 
the basic form of SOR-SVD is computed as follows: using a random number generator, we form a real matrix ${\bf \Omega} \in \mathbb R^{n \times \ell}$ with entries being independent, identically distributed (i.i.d.) Gaussian random variables of zero mean and unit variance. We then compute the matrix product: 
\begin{equation}
{\bf T}_1 = {\bf A}{\bf \Omega},
\label{eqT1}
\end{equation}
where ${\bf T}_1 \in \mathbb R^{m \times \ell}$ is formed by linear combinations of columns of $\bf A$ by the random Gaussian matrix. ${\bf T}_1$ is nothing but a projection onto the subspace spanned by columns of ${\bf A}$.
Having ${\bf T}_1$, we form the matrix ${\bf T}_2 \in \mathbb R^{n \times \ell}$: 
\begin{equation}
{\bf T}_2 = {\bf A}^T{\bf T}_1,
\label{eqAT1}
\end{equation}
where ${\bf T}_2$ is constructed by linear combinations of rows of $\bf A$ by ${\bf T}_1$. ${\bf T}_2$ is nothing but a projection onto the subspace spanned by rows of ${\bf A}$. Using the QR decomposition algorithm, we factor ${\bf T}_1$ and ${\bf T}_2$ such that: 
\begin{equation}
{\bf T}_1 = {\bf Q}_1{\bf R}_1  \quad \text{and} \quad {\bf T}_2 = 
{\bf Q}_2{\bf R}_2,
\end{equation}
where ${\bf Q}_1$ and ${\bf Q}_2$ are approximate bases for $\mathcal{R}({\bf A})$, 
and $\mathcal{R}({\bf A}^T)$, respectively.
We now form a matrix ${\bf M} \in \mathbb R^{\ell \times \ell}$ by compression of $\bf A$ through left and right multiplications by orthonormal bases:
\begin{equation}
{\bf M}={\bf Q}_1^T{\bf A}{\bf Q}_2,
\label{eqM}
\end{equation}

We then compute the rank-$k$ truncated SVD of ${\bf M}$:\looseness-1
\begin{equation}
{\bf M}_k = \widetilde{\bf U}_k \widetilde{\bf \Sigma}_k \widetilde{\bf V}_k,
\end{equation}
 Finally, we form the SOR-SVD-based low-rank approximation of $\bf A$: 
\begin{equation}
\hat{\bf A}_\text{SOR}= ({\bf Q}_1 \widetilde{\bf U}_k)\widetilde{\bf \Sigma}_k({\bf Q}_2  \widetilde{\bf V}_k)^T,
\end{equation}
where ${\bf Q}_1 \widetilde{\bf U}_k \in \mathbb R^{m \times k}$ and ${\bf Q}_2  \widetilde{\bf V}_k \in \mathbb R^{n \times k}$ are approximations to the $k$ leading left and right singular vectors of $\bf A$, respectively, and 
$\widetilde{\bf \Sigma}_k$ contains an approximation to the $k$ leading 
singular values. The algorithm is presented in Alg. \ref{Alg_first}. 

\begin{algorithm}
\caption{Subspace-Orbit Randomized SVD (SOR-SVD)}
\renewcommand{\algorithmicrequire}{\textbf{Input:}}
\begin{algorithmic}[1]
\REQUIRE ~~ 
 Matrix $\ {\bf A} \in \mathbb R^{m \times n}$, integer $\ell$,
\renewcommand{\algorithmicrequire}{\textbf{Output:}}
\REQUIRE ~~ A rank-$k$ approximation.
  \STATE Draw a standard Gaussian matrix ${\bf \Omega}\in \mathbb R^{n\times \ell}$;
  \STATE Compute ${\bf T}_1 = {\bf A}{\bf \Omega}$; \\
  \STATE Compute ${\bf T}_2 = {\bf A}^T{\bf T}_1$;
  \STATE Compute QR decompositions ${\bf T}_1 = {\bf Q}_1{\bf R}_1$, ${\bf T}_2 = {\bf Q}_2{\bf R}_2$; 
  \STATE Compute ${\bf M}={\bf Q}_1^T{\bf A}{\bf Q}_2$;
  \STATE Compute the rank-$k$ truncated SVD \\ ${\bf M}_k = \widetilde{\bf U}_k \widetilde{\bf \Sigma}_k \widetilde{\bf V}_k$;
  \STATE Form the SOR-SVD-based low-rank approximation of $\bf A$: $\hat{\bf A}_\text{SOR}= ({\bf Q}_1 \widetilde{\bf U}_k)\widetilde{\bf \Sigma}_k({\bf Q}_2  \widetilde{\bf V}_k)^T$.
\end{algorithmic}\label{Alg_first}
\end{algorithm}

The SOR-SVD requires three passes through data, for matrices stored-out-of-core, 
but it can be modified to revisit the data only once. To this end, in a similar manner as in Alg. \ref{Alg2}, the compressed matrix $\bf M$ in equation \eqref{eqM} can be computed by making use of currently available matrices as follows: both sides of the currently unknown ${\bf M}$ are postmultiplied by ${\bf Q}_2^T{\bf \Omega}$:
\begin{equation}
{\bf M}{\bf Q}_2^T{\bf \Omega} = {\bf Q}_1^T{\bf A}{\bf Q}_2{\bf Q}_2^T{\bf \Omega}.
\end{equation}

Having defined ${\bf A}\approx {\bf A}{\bf Q}_2{\bf Q}_2^T$ and ${\bf T}_1 = {\bf A}{\bf \Omega}$, an approximation to ${\bf M}$ can be obtained by
\begin{equation}
{\bf M}_\text{approx} = {\bf Q}_1^T{\bf T}_1({\bf Q}_2^T{\bf \Omega})^\dagger.
\end{equation} 

\begin{remark}
In the SOR-SVD algorithm, projecting $\bf A$ onto a
subspace spanned by its rows using a matrix containing random linear combinations 
of its columns, i.e., equation \eqref{eqAT1}, significantly improves the quality 
of the approximate basis ${\bf Q}_2$, compared to that of the TSR-SVD. This results 
in \textit{i)} a good approximation ${\bf M}_\text{approx}$ to $\bf M$, and \textit{ii)} tighter bounds for the singular values.
\end{remark}

The two key differences between our proposed SOR-SVD algorithm and the TSR-SVD algorithm are \textit{i)} to use a sketch of the input matrix in order to project it onto its row space, rather than using a random matrix, and \textit{ii)} to apply truncated SVD on the reduced-size matrix.

The SOR-SVD algorithm may be sufficiently accurate for matrices whose singular values display some decay, but, however, in applications where the data matrix has a slowly decaying singular values, it may produce poor singular vectors and singular values. Thus, we incorporate $q$ steps of a power iteration \cite{Rokhlin09,HMT2009}, to improve the accuracy of the algorithm in these circumstances. Given the matrix ${\bf A}$, and integers $k\le \ell< n$ and $q$, the resulting algorithm is described in Alg. \ref{Alg3}.\looseness-1

\begin{algorithm}
\caption{The SOR-SVD algorithm with Power Method}
\renewcommand{\algorithmicrequire}{\textbf{Input:}}
\begin{algorithmic}[1]
\REQUIRE ~~ 
 Matrix $\ {\bf A} \in \mathbb R^{m \times n}$,
integers $\ell$ and $q$,
\renewcommand{\algorithmicrequire}{\textbf{Output:}}
\REQUIRE ~~ A rank-$k$ approximation.
  \STATE Draw a standard Gaussian matrix ${\bf T}_2 \in \mathbb R^{n \times \ell}$;
  \FOR{$i=$ 1: $q+1$}
   \STATE Compute ${\bf T}_1 = {\bf A}{\bf T}_2$; \\
   \STATE Compute ${\bf T}_2 = {\bf A}^T{\bf T}_1$;
  \ENDFOR \\
  \STATE Compute QR decompositions ${\bf T}_1 = {\bf Q}_1{\bf R}_1$, 
  ${\bf T}_2 = {\bf Q}_2{\bf R}_2$; 
  \STATE Compute ${\bf M}={\bf Q}_1^T{\bf A}{\bf Q}_2$ or ${\bf M}_\text{approx} = {\bf Q}_1^T{\bf T}_1({\bf Q}_2^T{\bf T}_2)^
  \dagger $;
  \STATE Compute the rank-$k$ truncated SVD \\ ${\bf M}_k = \widetilde{\bf U}_k \widetilde{\bf \Sigma}_k \widetilde{\bf V}_k$ or ${\bf M}_\text{approx-k} = \widetilde{\bf U}_k \widetilde{\bf \Sigma}_k \widetilde{\bf V}_k$;
  \STATE Form the SOR-SVD-based low-rank approximation of $\bf A$: $\hat{\bf A}_\text{SOR}= ({\bf Q}_1 \widetilde{\bf U}_k)\widetilde{\bf \Sigma}_k({\bf Q}_2  \widetilde{\bf V}_k)^T$.
\end{algorithmic}\label{Alg3}
\end{algorithm}

Note that to compute the SOR-SVD when the power method is employed, a non-updated ${\bf T}_2$ must be used to form ${\bf M}_\text{approx}$.

\section{Analysis of SOR-SVD}
\label{secAnalysis}
In this section, we provide a detailed analysis of the performance of the SOR-SVD algorithms, the basic version Alg. \ref{Alg_first} as well as the one in Alg. \ref{Alg3}. In particular, we develop lower bounds for singular values, and upper bounds for rank-$k$ approximations in terms of both the spectral norm and the Frobenius norm.

While the TSR-SVD algorithm, described in Alg. \ref{Alg2}, results in the approximation ${\bf A} \approx {\bf Q}_1{\bf M}{\bf Q}_2^T$, where ${\bf M}$ is defined in \eqref{eqM}, the SOR-SVD procedure results in ${\bf A} \approx {\bf Q}_1{\bf M}_k{\bf Q}_2^T$, where ${\bf M}_k$ is the rank-$k$ truncated SVD of ${\bf M}$. The following theorem is a generalization of Theorem \ref{ThrEckYou} for the SOR-SVD.  

\begin{theorem}
Let ${\bf Q}_1 \in \mathbb R^{m \times \ell}$ and ${\bf Q}_2 \in \mathbb 
R^{n \times \ell}$ be matrices with orthonormal columns, and $1 \le k \le \ell$. 
Let furthermore ${\bf M}_k$ be the rank-$k$ truncated SVD of ${\bf Q}_1^T{\bf A}{\bf Q}_2$. Then ${\bf M}_k$ is an optimal solution to the following convex program:
\begin{equation}
\begin{aligned}
\underset{{\bf M} \in \mathbb R^{l \times l}, \text{rank}({\bf M})\le k}{\text{minimize}}
\|{\bf A} - {\bf Q}_1{\bf M}{\bf Q}_2^T\|_F = 
\|{\bf A} - {\bf Q}_1{\bf M}_k{\bf Q}_2^T\|_F,
\end{aligned}
\label{equThrFirst1}
\end{equation}
and  
\begin{equation}
\begin{aligned}
\|{\bf A} - {\bf Q}_1{\bf M}_k{\bf Q}_2^T\|_F \le &
\|{\bf A}_0\|_F +\|{\bf A}_k - {\bf Q}_1{\bf Q}_1^T{\bf A}_k\|_F \\
& + \|{\bf A}_k - {\bf A}_k{\bf Q}_2{\bf Q}_2^T\|_F,
\end{aligned}
\label{equThrFirst2}
\end{equation}
and we also have 
\begin{equation}
\begin{aligned}
\|{\bf A} - {\bf Q}_1{\bf M}_k{\bf Q}_2^T\|_2 \le &
\|{\bf A}_0\|_2 +\|{\bf A}_k - {\bf Q}_1{\bf Q}_1^T{\bf A}_k\|_F \\
& + \|{\bf A}_k - {\bf A}_k{\bf Q}_2{\bf Q}_2^T\|_F.
\end{aligned}
\label{equThrFirst3}
\end{equation}
\label{ThrFirst}
\end{theorem} 

\textit{Proof.}
We first rewrite the term in the left-hand side of \eqref{equThrFirst1}:
\begin{equation}
\begin{aligned}
& \|{\bf A} - {\bf Q}_1{\bf M}{\bf Q}_2^T\|_F^2 = \|{\bf A} - {\bf Q}_1{\bf Q}_1^T{\bf A}{\bf Q}_2{\bf Q}_2^T \\ & + {\bf Q}_1{\bf Q}_1^T{\bf A}{\bf Q}_2{\bf Q}_2^T 
- {\bf Q}_1{\bf M}{\bf Q}_2^T\|_F^2\\
& = \|{\bf A} - {\bf Q}_1{\bf Q}_1^T{\bf A}{\bf Q}_2{\bf Q}_2^T\|_F^2 + 
\|{\bf Q}_1^T{\bf A}{\bf Q}_2-{\bf M}\|_F^2.
\end{aligned}
\end{equation}
By Theorem \ref{ThrEckYou}, the result in \eqref{equThrFirst1} immediately follows.

According to Theorem \ref{ThrEckYou}, ${\bf A}_k$ is the best low-rank approximation to ${\bf A}$, whereas according to Theorem \ref{ThrFirst}, ${\bf Q}_1{\bf M}_k{\bf Q}_2^T$ is the best restricted (within a subspace) low-rank approximation to ${\bf A}$, with respect to the Frobenius norm. This leads to the following result 
\begin{equation}
\begin{aligned}
\|{\bf A} - {\bf A}_k\|_F \le & \|{\bf A} - {\bf Q}_1{\bf M}_k{\bf Q}_2^T\|_F\\
 \le & \|{\bf A} - {\bf Q}_1{\bf Q}_1^T{\bf A}_k{\bf Q}_2{\bf Q}_2^T\|_F.
\end{aligned} 
\label{equINEQ}
\end{equation}

The second relation holds because ${\bf Q}_1{\bf Q}_1^T{\bf A}_k{\bf Q}_2{\bf Q}_2^T$ is an undistinguished restricted Frobenius norm approximation to ${\bf A}$.

To prove \eqref{equThrFirst2}, we calculate
\begin{equation}
\begin{aligned}
& \|{\bf A} - {\bf Q}_1{\bf Q}_1^T{\bf A}_k{\bf Q}_2{\bf Q}_2^T\|_F^2 \\ 
& = \text{trace}\big(({\bf A} - {\bf Q}_1{\bf Q}_1^T{\bf A}_k{\bf Q}_2{\bf Q}_2^T)^T
({\bf A} - {\bf Q}_1{\bf Q}_1^T{\bf A}_k{\bf Q}_2{\bf Q}_2^T)\big) \\
& =  \text{trace}\big(({\bf A} - {\bf A}_k{\bf Q}_2{\bf Q}_2^T + {\bf A}_k{\bf Q}_2{\bf Q}_2^T - {\bf Q}_1{\bf Q}_1^T{\bf A}_k{\bf Q}_2{\bf Q}_2^T)^T \\
& \quad \times({\bf A} - {\bf A}_k{\bf Q}_2{\bf Q}_2^T + {\bf A}_k{\bf Q}_2{\bf Q}_2^T - {\bf Q}_1{\bf Q}_1^T{\bf A}_k{\bf Q}_2{\bf Q}_2^T)\big) \\
& = \|{\bf A} - {\bf A}_k{\bf Q}_2{\bf Q}_2^T\|_F^2 + 
\|{\bf A}_k{\bf Q}_2{\bf Q}_2^T -{\bf Q}_1{\bf Q}_1^T{\bf A}_k{\bf Q}_2{\bf Q}_2^T\|_F^2 + \\
& \quad 2\text{trace}\big( 
({\bf A} - {\bf A}_k{\bf Q}_2{\bf Q}_2^T)^T({\bf A}_k{\bf Q}_2{\bf Q}_2^T -{\bf Q}_1{\bf Q}_1^T{\bf A}_k{\bf Q}_2{\bf Q}_2^T)\big)\\
& = \|{\bf A} - {\bf A}_k{\bf Q}_2{\bf Q}_2^T\|_F^2 + 
\|{\bf A}_k -{\bf Q}_1{\bf Q}_1^T{\bf A}_k\|_F^2\\
& \quad + 2\text{trace}\big( 
({\bf I} - {\bf Q}_2{\bf Q}_2^T)\underbrace{{\bf A}_k{\bf Q}_2{\bf Q}_2^T({\bf A} - {\bf A}_k{\bf Q}_2{\bf Q}_2^T)^T}_{0}\big).
\end{aligned}
\label{equTRACE}
\end{equation}

Combining the last relation in \eqref{equTRACE} with equation \eqref{equINEQ}, gives
\begin{equation}
\begin{aligned}
\|{\bf A} - {\bf Q}_1{\bf M}_k{\bf Q}_2^T\|_F 
 \le & \|{\bf A} - {\bf A}_k{\bf Q}_2{\bf Q}_2^T\|_F \\ & + \|{\bf A}_k-{\bf Q}_1{\bf Q}_1^T{\bf A}_k\|_F.
\end{aligned} 
\label{}
\end{equation}

Writing ${\bf A} = {\bf A}_k + {\bf A}_0$, and applying the triangle inequality gives equation \eqref{equThrFirst2}.

To prove \eqref{equThrFirst3}, we observe that
\begin{equation}
\begin{aligned}
\|{\bf A} - {\bf Q}_1{\bf M}{\bf Q}_2^T\|_2 \le &
\|{\bf A} - {\bf Q}_1{\bf M}_k{\bf Q}_2^T\|_2 \\
\le & \|{\bf A} - {\bf Q}_1{\bf Q}_1^T{\bf A}_k{\bf Q}_2{\bf Q}_2^T\|_2\\
= & \|{\bf A}_0 + {\bf A}_k-{\bf Q}_1{\bf Q}_1^T{\bf A}_k{\bf Q}_2{\bf Q}_2^T\|_2\\ 
\le & \|{\bf A}_0\|_2 + \|{\bf A}_k - {\bf Q}_1{\bf Q}_1^T{\bf A}_k{\bf Q}_2{\bf Q}_2^T\|_2.
\end{aligned}
\label{equUnkn1}
\end{equation}

We can write the second term of the last equation as 
\begin{equation}
\begin{aligned}
& \|{\bf A}_k - {\bf Q}_1{\bf Q}_1^T{\bf A}_k{\bf Q}_2{\bf Q}_2^T\|_2\\
& = \|{\bf A}_k - {\bf Q}_1{\bf Q}_1^T{\bf A}_k + {\bf Q}_1{\bf Q}_1^T{\bf A}_k 
- {\bf Q}_1{\bf Q}_1^T{\bf A}_k{\bf Q}_2{\bf Q}_2^T\|_2 \\
& \le \|{\bf A}_k - {\bf Q}_1{\bf Q}_1^T{\bf A}_k\|_2 +  \|{\bf Q}_1{\bf Q}_1^T({\bf A}_k -{\bf A}_k{\bf Q}_2{\bf Q}_2^T)\|_2 \\
& \le \|{\bf A}_k - {\bf Q}_1{\bf Q}_1^T{\bf A}_k\|_F +  \|{\bf A}_k -
{\bf A}_k{\bf Q}_2{\bf Q}_2^T|_F.
\end{aligned}
\label{equUnkn2}
\end{equation}

Plugging this into equation \eqref{equUnkn1} yields \eqref{equThrFirst3}. 

The last relation in \eqref{equUnkn2} holds due to the unitary invariance 
property of the $\ell_2$-norm, and furthermore, to the relation that for 
any matrix ${\bf \Pi}$, $\|{\bf \Pi}\|_2 \le \|{\bf \Pi}\|_F$. \QEDB

\subsection{Deterministic Error Bounds}
\label{secDetEB}

In this section, we make use of techniques from linear algebra to give generic error bounds which depend on the interaction between the standard Gaussian matrix $\bf \Omega$ and the right singular vectors of the data matrix $\bf A$. 

To derive lower bounds on approximated singular values, we begin by stating two key results that are used in the analysis later on.

\begin{theorem}(Thompson \cite{Thompson72})
Suppose that the matrix $\bf A$ has singular values as defined in \eqref{eqSVD}. Let ${\bf M} \in \mathbb R^{\ell \times \ell}$ be a submatrix of $\bf A$. 
Then for $j = 1,..., \ell,$ we have 
\begin{equation}
\sigma_j \ge \sigma_j(\bf M).
\label{eqThrThomp}
\end{equation}
\end{theorem}

The relation in \eqref{eqThrThomp} can be easily proven by allowing $\bf M$ to be ${\bf M}={\bf Q}_1^T{\bf A}{\bf Q}_2$, where ${\bf Q}_1 \in \mathbb R^{m \times \ell}$ and ${\bf Q}_2 \in \mathbb R^{n \times \ell}$ are orthonormal matrices.  

\begin{remark}
Since the matrices ${\bf M}$ and ${\bf M}_k$ have the same singular values  
$\sigma_j$, for $j = 1,..., k$, and moreover, singular values of
${\bf Q}_1^T{\bf A}{\bf Q}_2$ and ${\bf Q}_1{\bf Q}_1^T{\bf A}{\bf Q}_2
{\bf Q}_2^T$ coincide \cite{Thompson72}, we have 
\begin{equation}
\sigma_j \ge \sigma_j({\bf Q}_1{\bf M}_k{\bf Q}_2^T) = \sigma_j({\bf Q}_1{\bf Q}_1^T{\bf A}{\bf Q}_2{\bf Q}_2^T).
\label{eqRemark}
\end{equation} 
\label{RemarkThomp}
\end{remark}
 
\begin{lemma}(Gu \cite{Gu2015}). Given any matrix 
${\bf B} \in \mathbb R^{m \times n}$, and a  full column rank matrix 
${\bf \Omega}\in \mathbb R^{n \times \ell}$, assume that ${\bf X}\in 
\mathbb R^{\ell \times \ell}$ is a non-singular matrix. Let ${\bf B}{\bf \Omega} 
= {\bf Q}_\text{b}{\bf R}_\text{b}$ and ${\bf B}{\bf \Omega X} = {\bf Q}_\text{x}{\bf R}_\text{x}$ be QR decompositions of the matrix products. Then 
\begin{equation}
{\bf Q}_\text{b}{\bf Q}_\text{b}^T = {\bf Q}_\text{x}{\bf Q}_\text{x}^T. 
\end{equation}
\label{Lemma1}
\end{lemma} 

For the matrix ${\bf T}_2$, equation \eqref{eqAT1}, and by the SVD of 
$\bf A$, given in \eqref{eqSVD}, we have
\begin{equation}
{\bf T}_2  = {\bf A}^T{\bf T}_1 = {\bf A}^T{\bf A \Omega} = 
{\bf V \Sigma}^2{\bf V}^T \bf \Omega.
\end{equation}

By partitioning ${\bf \Sigma}^2$, we have
\begin{equation}
\begin{aligned}
{\bf T}_2 
 ={\bf V}\begin{bmatrix}
       {\bf \Sigma}_1^{2} & {\bf 0} & {\bf 0} \\
       {\bf 0} & {\bf \Sigma}_2^{2} & {\bf 0}\\
       {\bf 0} & {\bf 0} & {\bf \Sigma}_3^{2}
  \end{bmatrix}
  \begin{bmatrix}
       {\bf V}_1^T{\bf \Omega} \\
       \\
       {\bf V}_2^T {\bf \Omega}
       \end{bmatrix} = {\bf Q}_2{\bf R}_2,
\end{aligned}
\end{equation}
where ${\bf \Sigma}_1 \in \mathbb R^{k \times k}$, ${\bf \Sigma}_2 \in 
\mathbb R^{\ell-p-k \times \ell-p-k}$, ${\bf \Sigma}_3 \in \mathbb 
R^{n-\ell+p \times n-\ell+p}$.

We define ${\bf \Omega}_1 \in \mathbb R^{\ell-p \times \ell}$, 
and ${\bf \Omega}_2 \in \mathbb R^{n-\ell+p \times \ell}$ as follows: 
\begin{equation}
{\bf \Omega}_1  \triangleq {\bf V}_1^T{\bf \Omega}, \quad \text{and} \quad
{\bf \Omega}_2  \triangleq {\bf V}_2^T{\bf \Omega}.
\label{eqOmegas}
\end{equation}

We assume that ${\bf \Omega}_1$ is full row rank and its pseudo-inverse satisfies
\begin{equation}
{\bf \Omega}_1{\bf \Omega}_1^\dagger  =\bf I.
\end{equation}

To understand the behavior of singular values and the low-rank approximation, we choose a matrix ${\bf X} \in \mathbb R^{\ell \times \ell}$, which orients the first $k$ columns of ${\bf T}_2{\bf X}$ in the directions of the $k$ leading singular vectors in $\bf V$. Thus we choose
\begin{equation}
{\bf X} = \bigg[
{\bf \Omega}_1^\dagger
   \begin{pmatrix} 
        {\bf \Sigma}_1^{2} & {\bf 0}\\
        {\bf 0} & {\bf \Sigma}_2^{2}
   \end{pmatrix}^{-1}, \tilde{\bf X}\bigg],
   \label{eqX_1}
\end{equation}
where the $\tilde{\bf X}\in \mathbb R^{\ell \times p}$ is chosen such that ${\bf X}\in \mathbb R^{\ell \times \ell}$ is non-singular and ${\bf \Omega}_1\tilde{\bf X} = \bf 0$. Now we can write
\begin{equation}
\begin{aligned}
{\bf T}_2{\bf X} = {\bf A}^T{\bf A \Omega X} = 
& {\bf V}\begin{bmatrix}
   \begin{pmatrix}
        {\bf \Sigma}_1^2 & {\bf 0}\\
        {\bf 0} & {\bf \Sigma}_2^2
   \end{pmatrix}{\bf \Omega}_1\\
   {\bf \Sigma}_3^2 {\bf \Omega}_2
   \end{bmatrix}{\bf X}\\
= & {\bf V} \begin{bmatrix}
       {\bf I} & {\bf 0} & {\bf 0} \\
       {\bf 0} & {\bf I} & {\bf 0}\\
       {\bf W}_1 & {\bf W}_2 & {\bf W}_3
  \end{bmatrix},
\end{aligned} 
\label{eqAtT1X}
\end{equation}
where ${\bf W}_1 = {\bf \Sigma}_3^2{\bf \Omega}_2{\bf \Omega}_1^\dagger {\bf \Sigma}_1^{-2} \in \mathbb R^{n-\ell+p \times k}$, ${\bf W}_2 = {\bf \Sigma}_3^2{\bf \Omega}_2{\bf \Omega}_1^\dagger {\bf \Sigma}_2^{-2} \in \mathbb R^{n-\ell+p \times \ell-p-k}$, and ${\bf W}_3 = {\bf \Sigma}_3^2{\bf \Omega}_2 \tilde{\bf X} \in \mathbb R^{n-\ell+p \times p}$.

Now by a QR decomposition of equation \eqref{eqAtT1X}, with partitioned 
matrices, we have
\begin{equation}
\begin{aligned}
& {\bf V} \begin{bmatrix}
       {\bf I} & {\bf 0} & {\bf 0} \\
       {\bf 0} & {\bf I} & {\bf 0}\\
       {\bf W}_1 & {\bf W}_2 & {\bf W}_3
  \end{bmatrix} =
  \tilde{\bf Q}\tilde{\bf R}\\
  & = [\tilde{\bf Q}_1 \quad \tilde{\bf Q}_2 \quad \tilde{\bf Q}_3] 
  \begin{bmatrix}
       \tilde{\bf R}_{11} & \tilde{\bf R}_{12} & \tilde{\bf R}_{13} \\
       {\bf 0} & \tilde{\bf R}_{22} & \tilde{\bf R}_{23}\\
       {\bf 0} & {\bf 0} & \tilde{\bf R}_{33}
  \end{bmatrix}\\
  & = \begin{bmatrix}
       (\tilde{\bf Q}_1\tilde{\bf R}_{11})^T \\
       (\tilde{\bf Q}_1\tilde{\bf R}_{12} + \tilde{\bf Q}_2\tilde{\bf R}_{22})^T \\
       (\tilde{\bf Q}_1\tilde{\bf R}_{13} + \tilde{\bf Q}_2\tilde{\bf R}_{23} + \tilde{\bf Q}_3\tilde{\bf R}_{33})^T 
  \end{bmatrix}^T.
\end{aligned} 
\label{eqRelations}
\end{equation}

We use this representation to develop lower bounds on singular values and upper bounds for the rank-$k$ approximation of the SOR-SVD algorithm.

\begin{lemma}
Let ${\bf W}_1$ be defined as in \eqref{eqAtT1X}, and assume that ${\bf \Omega}_1$ is full row rank. Then for $\hat{\bf A}_\text{SOR}$ we have
\begin{equation}
\sigma_k(\hat{\bf A}_\text{SOR}) \ge \frac{\sigma_k}{\sqrt{1 + {\|{{\bf W}_1}\|_2^2}}}.
\label{eqLemmW}
\end{equation}
\label{Lemm_sv_1st}
\end{lemma}

\textit{Proof.} 
According to Lemma \ref{Lemma1}, we have 
\begin{equation}
{\bf Q}_2{\bf Q}_2^T  =  \tilde{\bf Q}\tilde{\bf Q}^T.
\end{equation}

Thus
\begin{equation}
\begin{aligned}
{\bf Q}_1{\bf Q}_1^T{\bf A}{\bf Q}_2{\bf Q}_2^T & =  {\bf Q}_1{\bf Q}_1^T{\bf A}\tilde{\bf Q}\tilde{\bf Q}^T.
\end{aligned}
\end{equation}

We now write 
\begin{equation}
\begin{aligned}
{\bf A}\tilde{\bf Q}\tilde{\bf Q}^T& = {\bf A}(\tilde{\bf Q}_1 \quad \tilde{\bf Q}_2 \quad \tilde{\bf Q}_3)\tilde{\bf Q}^T\\
& = {\bf U}\underbrace{\begin{bmatrix}
       {\bf \Sigma}_1 & {\bf 0} & {\bf 0} \\
       {\bf 0} & {\bf \Sigma}_2 & {\bf 0}\\
       {\bf 0} & {\bf 0} & {\bf \Sigma}_3
  \end{bmatrix}{\bf V}^T[\tilde{\bf Q}_1 \quad \tilde{\bf Q}_2 \quad \tilde{\bf Q}_3]}_{\bf S}
  \tilde{\bf Q}.
\end{aligned}
\end{equation}

We now write $\bf S$ as:
\begin{equation}
\renewcommand\arraystretch{1.4}
\mleft[
\begin{array}{c|c}
  ({\bf \Sigma}_1 \quad {\bf 0} \quad {\bf 0}){\bf V}^T\tilde{\bf Q}_1 & 
  ({\bf \Sigma}_1 \quad {\bf 0} \quad {\bf 0}){\bf V}^T(\tilde{\bf Q}_2 \mspace{10mu} \tilde{\bf Q}_3)\\
  \hline
  \begin{pmatrix}
        {\bf 0} & {\bf 0}\\
        {\bf \Sigma}_1 & {\bf 0}\\
        {\bf 0} & {\bf \Sigma}_2
   \end{pmatrix}^T{\bf V}^T\tilde{\bf Q}_1 &  
 \begin{pmatrix}
        {\bf 0} & {\bf 0}\\
        {\bf \Sigma}_1 & {\bf 0}\\
        {\bf 0} & {\bf \Sigma}_2        
 \end{pmatrix}^T{\bf V}^T
   (\tilde{\bf Q}_2 \mspace{10mu} \tilde{\bf Q}_3)
\end{array}
\mright].
\notag
\end{equation}

We observe that the matrix $({\bf \Sigma}_1 \quad {\bf 0} \quad {\bf 0}){\bf V}^T\tilde{\bf Q}_1$ is a submatrix of ${\bf Q}_1{\bf Q}_1^T{\bf A}{\bf Q}_2{\bf Q}_2^T$. By relations in equation \eqref{eqRelations}, we have
\begin{equation}
{\bf V} \begin{bmatrix}
       {\bf I} \\
       {\bf 0} \\
       {\bf W}_1
  \end{bmatrix} = \tilde{\bf Q}_1\tilde{\bf R}_{11}. 
  \label{eqQ1tilde}
\end{equation}
and accordingly
\begin{equation}
\|\tilde{\bf R}_{11}\|_2=\sqrt{1+\|{\bf W}_1\|_2^2}.
\label{eqR_11}
\end{equation}

By substituting $\tilde{\bf Q}_1$, from \eqref{eqQ1tilde}, into $({\bf \Sigma}_1 \quad {\bf 0} \quad {\bf 0}){\bf V}^T\tilde{\bf Q}_1$,
\begin{equation}
\begin{aligned}
&({\bf \Sigma}_1 \quad {\bf 0} \quad {\bf 0}){\bf V}^T\tilde{\bf Q}_1 \\
& = ({\bf \Sigma}_1 \quad {\bf 0} \quad {\bf 0}){\bf V}^T{\bf V} 
   \begin{bmatrix}
       {\bf I} \\
       {\bf 0} \\
       {\bf W}_1
   \end{bmatrix} \tilde{\bf R}_{11}^{-1} \\ 
& = {\bf \Sigma}_1\tilde{\bf R}_{11}^{-1},
\end{aligned}
\end{equation}
and by Remark \ref{RemarkThomp} it follows that
\begin{equation}
\begin{aligned}
&\sigma_j(\hat{\bf A}_\text{SOR}) = 
\sigma_j({\bf Q}_1{\bf Q}_1^T{\bf A}{\bf Q}_2{\bf Q}_2^T)\ge \sigma_j({\bf \Sigma}_1\tilde{\bf R}_{11}^{-1}).
\end{aligned}
\label{eqResultRem}
\end{equation}

On the other hand, we also have
\begin{equation}
\sigma_j = \sigma_j ({\bf \Sigma}_1\tilde{\bf R}_{11}^{-1}\tilde{\bf R}_{11})
\le \sigma_j ({\bf \Sigma}_1\tilde{\bf R}_{11}^{-1})\|\tilde{\bf R}_{11}\|_2.
\end{equation}

Plugging the last relation into \eqref{eqResultRem} and using \eqref{eqR_11}, we obtain \eqref{eqLemmW}. \QEDB 

We are now prepared to present the result.

\begin{theorem}
Suppose that the matrix $\bf A$ has an SVD defined in \eqref{eqSVD}, $2\le p+k\le \ell$, and $\hat{\bf A}_\text{SOR}$ is computed through the basic form of SOR-SVD, Alg. \ref{Alg_first}. Assume furthermore that ${\bf \Omega}_1$ is full row rank, then for $j = 1,..., k,$ we have 
\begin{equation}
\sigma_j \ge \sigma_j(\hat{\bf A}_\text{SOR}) \ge \frac{\sigma_j}{\sqrt{1 + {\|{{\bf \Omega}_2}\|_2^2} \|{{\bf \Omega}_1}^\dagger\|_2^2
\Big(\frac{\sigma_{\ell-p+1}}{\sigma_j}\Big)^4}}.
\label{eqThr1_sv1}
\end{equation}
and when the power method is used, i.e., Alg. \ref{Alg3}, we have 
\begin{equation}
\sigma_j \ge \sigma_j(\hat{\bf A}_\text{SOR}) \ge \frac{\sigma_j}{\sqrt{1 + {\|{{\bf \Omega}_2}\|_2^2} \|{{\bf \Omega}_1}^\dagger\|_2^2
\Big(\frac{\sigma_{\ell-p+1}}{\sigma_j}\Big)^{4q+4}}}.
\label{eqThr1_sv2}
\end{equation}
\label{Thr1}
\end{theorem}

\textit{Proof.}
To prove \eqref{eqThr1_sv1}, according to the definition of ${\bf W}_1$ in 
equation \eqref{eqAtT1X}, we obtain
\begin{equation}
\begin{aligned}
\|{\bf W}_1\|_2 \le \Big(\frac{\sigma_{\ell-p+1}}{\sigma_k}\Big)^2
\|{\bf \Omega}_2\|_2\|{\bf \Omega}_1^\dagger\|_2.
\end{aligned}
\end{equation}

Taking this together with equation \eqref{eqLemmW}, yields the result in \eqref{eqThr1_sv1} for $j=k$. 

To prove \eqref{eqThr1_sv2}, we observe that when the power method is used, 
the SOR-SVD computation starts off by setting ${\bf T}_2 = {\bf \Omega}$, 
and ${\bf T}_2$ is updated such that
\begin{equation}
{\bf T}_2 = ({\bf A}^T{\bf A})^q{\bf A}^T{\bf T}_1 = ({\bf A}^T{\bf A})^q{\bf A}^T{\bf A \Omega}.
\end{equation}

By writing the SVD of $\bf A$ \eqref{eqSVD},
\begin{equation}
\begin{aligned}
{\bf T}_2 & = {\bf V \Sigma}^{2q+2}{\bf V}^T \bf \Omega\\
& ={\bf V}\begin{bmatrix}
       {\bf \Sigma}_1^{2q+2} & {\bf 0} & {\bf 0} \\
       {\bf 0} & {\bf \Sigma}_2^{2q+2} & {\bf 0}\\
       {\bf 0} & {\bf 0} & {\bf \Sigma}_3^{2q+2}
  \end{bmatrix}
  \begin{bmatrix}
       {\bf V}_1^T{\bf \Omega} \\
       \\
       {\bf V}_2^T {\bf \Omega}
       \end{bmatrix} = {\bf Q}_2{\bf R}_2.
\end{aligned}
\end{equation}

Consequently, the matrix $\bf X$, defined in \eqref{eqX_1}, is now 
defined as follows:
\begin{equation}
{\bf X} = \bigg[
{\bf \Omega}_1^\dagger
   \begin{pmatrix}
        {\bf \Sigma}_1^{2q+2} & {\bf 0}\\
        {\bf 0} & {\bf \Sigma}_2^{2q+2}
   \end{pmatrix}^{-1}, \tilde{\bf X}\bigg].
\end{equation}

Forming ${\bf T}_2{\bf X}$ yields:
\begin{equation}
{\bf W}_1 = {\bf \Sigma}_3^{2q+2}{\bf \Omega}_2{\bf \Omega}_1^\dagger 
{\bf \Sigma}_1^{-(2q+2)}.
\label{eqW1_PM}
\end{equation}

Thus 
\begin{equation}
\begin{aligned}
\|{\bf W}_1\|_2 \le \Big(\frac{\sigma_{\ell-p+1}}{\sigma_k}\Big)^{2q+2}
\|{\bf \Omega}_2\|_2\|{\bf \Omega}_1^\dagger\|_2.
\end{aligned}
\end{equation}

Taking this together with equation \eqref{eqLemmW}, yields the result for $j=k$. 

To prove Theorem \ref{Thr1} for any $1 \le j < k$, since by Remark 
\ref{RemarkThomp}, $\sigma_j(\hat{\bf A}_\text{SOR}) = \sigma_j({\bf Q}_1{\bf Q}_1^T{\bf A}{\bf Q}_2{\bf Q}_2^T)$, it is only required to repeat all 
previous arguments for a rank $j$ truncated SVD. \QEDB

To derive lower bounds on the low-rank approximation, with respect to both the spectral norm and the Frobenius norm, provided equations \eqref{equThrFirst2}, \eqref{equThrFirst3}, it suffices to bound the right-hand sides of the equations, i.e., $\|{\bf A}_k - {\bf Q}_1{\bf Q}_1^T{\bf A}_k\|_F$ and $\|{\bf A}_k - {\bf A}_k{\bf Q}_2{\bf Q}_2^T\|_F$. To this end, we begin by stating a proposition from \cite{HMT2009}.

\begin{proposition}(Halko et al. \cite{HMT2009})
Suppose that for given matrices ${\bf N}_1$ and ${\bf N}_2$, $\mathcal{R}({\bf N}_1) \subset \mathcal{R}({\bf N}_2)$. Then for any matrix 
$\bf A$, it holds that
\begin{equation}
\begin{aligned}
& \|{\bf P}_{{\bf N}_1}{\bf A}\|_2 \le \|{\bf P}_{{\bf N}_2}{\bf A}\|_2, \\
& \|({\bf I} - {\bf P}_{{\bf N}_2}){\bf A}\|_2 \le \|({\bf I}-{\bf P}_{{\bf N}_1}){\bf A}\|_2,
\end{aligned}
\end{equation}
where ${\bf P}_{{\bf N}_1}$ and ${\bf P}_{{\bf N}_2}$ are orthogonal projections onto ${\bf N}_1$ and ${\bf N}_2$, respectively.
\label{PropHalkoP}
\end{proposition}

By combining Lemma \ref{Lemma1} and Proposition \ref{PropHalkoP}, it follows that
\begin{equation}
\begin{aligned}
\|{\bf A}_k({\bf I} - {\bf Q}_2{\bf Q}_2^T)\|_F 
& = \|{\bf A}_k({\bf I} - \tilde{\bf Q}\tilde{\bf Q}^T)\|_F \\
& \le \|{\bf A}_k({\bf I} - \tilde{\bf Q}_1 \tilde{\bf Q}_1^T)\|_F.
\end{aligned}
\label{eqcomL1P1_1}
\end{equation}

By the definition of $\tilde{\bf Q}_1$, equation \eqref{eqQ1tilde}, it follows 
\begin{equation}
\begin{aligned}
&{\bf I} - \tilde{\bf Q}_1\tilde{\bf Q}_1^T = \\
& {\bf V} 
   \begin{bmatrix}
       {\bf I} - \tilde{\bf W}^{-1} & {\bf 0} & -\tilde{\bf W}^{-1}{\bf W}_1^T \\
       {\bf 0} & -{\bf I} & {\bf 0}\\
       - {\bf W}_1\tilde{\bf W}^{-1} & {\bf 0} & 
       {\bf I} - {\bf W}_1\tilde{\bf W}^{-1}{\bf W}_1^T
   \end{bmatrix} {\bf V}^T,
\end{aligned}
\end{equation}
where ${\bf W}_1$ is defined in \eqref{eqAtT1X}, and 
$\tilde{\bf W}^{-1}$ is defined as follows:
\begin{equation}
\begin{aligned}
\tilde{\bf R}_{11}^{-1}\tilde{\bf R}_{11}^{-T} = 
  (\tilde{\bf R}_{11}^T\tilde{\bf R}_{11})^{-1} = (1+{\bf W}_1^T{\bf W}_1) ^{-1} 
  = \tilde{\bf W}^{-1}.
\end{aligned}
\notag
\end{equation}

We then have 
\begin{equation}
\begin{aligned}
& {\bf A}_k({\bf I} - \tilde{\bf Q}_1\tilde{\bf Q}_1^T) = {\bf U}[{\bf \Sigma}_1 \quad 
{\bf 0}\quad {\bf 0}]{\bf V}^T({\bf I} - \tilde{\bf Q}_1\tilde{\bf Q}_1^T)\\
& = \underbrace{
{\bf U}[{\bf \Sigma}_1({\bf I}-\tilde{\bf W}^{-1}) \quad {\bf 0}\quad -
{\bf \Sigma}_1\tilde{\bf W}^{-1}{\bf W}_1^T]{\bf V}^T}_{\bf N}.
\end{aligned}
\end{equation}

By the definition of the Frobenius norm, it follows that  
\begin{equation}
\begin{aligned}
\|{\bf A}_k({\bf I} - & \tilde{\bf Q}_1\tilde{\bf Q}_1^T)\|_F = 
\sqrt{\text{trace}({{\bf NN}^T})}\\
= & \sqrt{\text{trace}({\bf \Sigma}_1{\bf W}_1^T({\bf I} + {\bf W}_1{\bf W}_1^T)^{-1}{\bf W}_1{\bf \Sigma}_1^T)}\\
\le & \sqrt{\text{trace}([\|{\bf W}_1{\bf \Sigma}_1\|_2^2{\bf I} + {\bf \Sigma}_1^{-2}]^{-1})}\\
= & \|{\bf W}_1{\bf \Sigma}_1\|_2 \times \\
& \sqrt{\text{trace}({\bf \Sigma}_1[\|{\bf W}_1{\bf \Sigma}_1\|_2^2{\bf I} + {\bf \Sigma}_1^{-2}]^{-1}){\bf \Sigma}_1^T)}\\
\le & \dfrac{\sqrt{k}\|{\bf W}_1{\bf \Sigma}_1\|_2 \sigma_1} {\sqrt{{\sigma_1}^2 + \|{\bf W}_1{\bf \Sigma}_1\|_2^2}}.
\end{aligned} 
\end{equation}

Plugging this into \eqref{eqcomL1P1_1}, we have 
\begin{equation}
\begin{aligned}
\|{\bf A}_k({\bf I} - {\bf Q}_2{\bf Q}_2^T)\|_F 
\le \dfrac{\sqrt{k}\|{\bf W}_1{\bf \Sigma}_1\|_2 \sigma_1} {\sqrt{{\sigma_1}^2 + \|{\bf W}_1{\bf \Sigma}_1\|_2^2}}.
\end{aligned}
\label{eqcomL1P1_WSigma}
\end{equation}
 
To bound $\|{\bf A}_k - {\bf Q}_1{\bf Q}_1^T{\bf A}_k\|_F$, we need to perform the same procedure described for ${\bf T}_2$, for the matrix product 
${\bf T}_1$ in equation \eqref{eqT1}. Thus, for ${\bf T}_1$, and by the SVD of 
$\bf A$ in \eqref{eqSVD}, we have
\begin{equation}
{\bf T}_1 = {\bf A}{\bf \Omega} = {\bf U \Sigma}{\bf V}^T \bf \Omega.
\label{eqT1_2}
\end{equation}

By partitioning ${\bf \Sigma}$, we obtain 
\begin{equation}
\begin{aligned}
{\bf T}_1
 ={\bf U}\begin{bmatrix}
       {\bf \Sigma}_1 & {\bf 0} & {\bf 0} \\
       {\bf 0} & {\bf \Sigma}_2 & {\bf 0}\\
       {\bf 0} & {\bf 0} & {\bf \Sigma}_3
  \end{bmatrix}
  \begin{bmatrix}
       {\bf V}_1^T{\bf \Omega} \\
       \\
       {\bf V}_2^T {\bf \Omega}
       \end{bmatrix} = {\bf Q}_1{\bf R}_1,
\end{aligned}
\label{eqT1_Q1}
\end{equation}

Having \eqref{eqOmegas}, we now choose a matrix ${\bf X} \in \mathbb R^{\ell \times \ell}$, which orients the first $k$ columns of ${\bf T}_1{\bf X}$ in the directions of the $k$ leading singular vectors in $\bf U$. Thus, we have
\begin{equation}
{\bf X} = \bigg[
{\bf \Omega}_1^\dagger
   \begin{pmatrix}
        {\bf \Sigma}_1 & {\bf 0}\\
        {\bf 0} & {\bf \Sigma}_2
   \end{pmatrix}^{-1}, \bar{\bf X}\bigg],
\end{equation}
where the $\bar{\bf X}\in \mathbb R^{\ell \times p}$ is chosen such that ${\bf X}\in 
\mathbb R^{\ell \times \ell}$ is non-singular and ${\bf \Omega}_1\bar{\bf X} = \bf 0$. 
We now write
\begin{equation}
\begin{aligned}
{\bf T}_1{\bf X}={\bf A \Omega X} = 
& {\bf U}\begin{bmatrix}
   \begin{pmatrix}
        {\bf \Sigma}_1 & {\bf 0}\\
        {\bf 0} & {\bf \Sigma}_2
   \end{pmatrix}{\bf \Omega}_1\\
   {\bf \Sigma}_3 {\bf \Omega}_2
   \end{bmatrix}{\bf X}\\
= & {\bf U} \begin{bmatrix}
       {\bf I} & {\bf 0} & {\bf 0} \\
       {\bf 0} & {\bf I} & {\bf 0}\\
       {\bf D}_1 & {\bf D}_2 & {\bf D}_3
  \end{bmatrix},    
\end{aligned} 
\label{eqD1_NoPM}
\end{equation}
where ${\bf D}_1 = {\bf \Sigma}_3{\bf \Omega}_2{\bf \Omega}_1^\dagger {\bf \Sigma}_1^{-1} \in \mathbb R^{n-\ell+p \times k}$, ${\bf D}_2 = {\bf \Sigma}_3{\bf \Omega}_2{\bf \Omega}_1^\dagger {\bf \Sigma}_2^{-1} \in \mathbb R^{n-\ell+p \times \ell-p-k}$, and ${\bf D}_3 = {\bf \Sigma}_3{\bf \Omega}_2 \bar{\bf X} \in \mathbb R^{n-\ell+p \times p}$.

We then write
\begin{equation}
\begin{aligned}
& {\bf U} \begin{bmatrix}
       {\bf I} & {\bf 0} & {\bf 0} \\
       {\bf 0} & {\bf I} & {\bf 0}\\
       {\bf D}_1 & {\bf D}_2 & {\bf D}_3
  \end{bmatrix} =
  \bar{\bf Q}\bar{\bf R}\\
  & = [\bar{\bf Q}_1 \quad \bar{\bf Q}_2 \quad \bar{\bf Q}_3] 
  \begin{bmatrix}
       \bar{\bf R}_{11} & \bar{\bf R}_{12} & \bar{\bf R}_{13} \\
       {\bf 0} & \bar{\bf R}_{22} & \bar{\bf R}_{23}\\
       {\bf 0} & {\bf 0} & \bar{\bf R}_{33} 
  \end{bmatrix}.
\end{aligned} 
\label{eqT1X_Qbar}
\end{equation}
and as a result, we have
\begin{equation}
\begin{aligned}
{\bf U} \begin{bmatrix}
       {\bf I} \\
       {\bf 0} \\
       {\bf D}_1
  \end{bmatrix} = \bar{\bf Q}_1\bar{\bf R}_{11}.
\end{aligned} 
\label{eqQ1bar}
\end{equation}

By Lemma \ref{Lemma1}, we have 
\begin{equation}
{\bf Q}_1{\bf Q}_1^T  =  \bar{\bf Q}\bar{\bf Q}^T.
\end{equation}
where ${\bf Q}_1$ and $\bar{\bf Q}$ are the Q-factors of the QR decompositions 
of ${\bf T}_1$ \eqref{eqT1_Q1}, and ${\bf T}_1{\bf X}$ \eqref{eqT1X_Qbar}, respectively. By combining Lemma \ref{Lemma1} and Proposition \ref{PropHalkoP}, 
it follows that
\begin{equation}
\begin{aligned}
\|({\bf I} - {\bf Q}_1{\bf Q}_1^T){\bf A}_k\|_F 
& = \|({\bf I} - \bar{\bf Q}\bar{\bf Q}^T){\bf A}_k\|_F \\
& \le \|({\bf I} - \bar{\bf Q}_1 \bar{\bf Q}_1^T){\bf A}_k\|_F.
\end{aligned}
\label{eqcomL1P1_second}
\end{equation}

By the definition of $\bar{\bf Q}_1$, equation \eqref{eqQ1bar}, it follows
\begin{equation}
\begin{aligned}
&{\bf I} - \bar{\bf Q}_1 \bar{\bf Q}_1^T = \\
& {\bf U} 
   \begin{bmatrix}
       {\bf I} - \bar{\bf D}^{-1} & {\bf 0} & -\bar{\bf D}^{-1}{\bf D}_1^T \\
       {\bf 0} & -{\bf I} & {\bf 0}\\
       - {\bf D}_1 \bar{\bf D}^{-1} & {\bf 0} & 
       {\bf I} - {\bf D}_1 \bar{\bf D}^{-1}{\bf D}_1^T
   \end{bmatrix} {\bf U}^T,
\end{aligned}
\end{equation}
where ${\bf D}_1$ is defined in \eqref{eqD1_NoPM}, and $\bar{\bf D}^{-1}$ is 
defined as follows:
\begin{equation}
\begin{aligned}
\bar{\bf R}_{11}^{-1} \bar{\bf R}_{11}^{-T} = 
  (\bar{\bf R}_{11}^T \bar{\bf R}_{11})^{-1} = (1+{\bf D}_1^T{\bf D}_1) ^{-1} 
  = \bar{\bf D}^{-1}.
\end{aligned}
\notag
\end{equation}

We then write 
\begin{equation}
\begin{aligned}
({\bf I} - \bar{\bf Q}_1\bar{\bf Q}_1^T){\bf A}_k & = ({\bf I} - \bar{\bf Q}_1\bar{\bf Q}_1^T){\bf U}\begin{bmatrix}
       {\bf \Sigma}_1  \\
       {\bf 0} \\
       {\bf 0}
   \end{bmatrix}{\bf V}^T\\
& = \underbrace{{\bf U} 
   \begin{bmatrix}
       ({\bf I} - \bar{\bf D}^{-1}){\bf \Sigma}_1  \\
       {\bf 0} \\
       - {\bf D}_1 \bar{\bf D}^{-1}{\bf \Sigma}_1
   \end{bmatrix} {\bf V}^T}_{\bf H}.
\end{aligned}
\end{equation}

By the definition of the Frobenius norm, it follows that 
\begin{equation}
\begin{aligned}
\|({\bf I} - & \bar{\bf Q}_1 \bar{\bf Q}_1^T){\bf A}_k\|_F = 
\sqrt{\text{trace}({{\bf H}^T{\bf H}})}\\
= & \sqrt{\text{trace}({\bf \Sigma}_1{\bf D}_1^T({\bf I} + {\bf D}_1{\bf D}_1^T)^{-1}{\bf D}_1{\bf \Sigma}_1^T)}\\
\le & \sqrt{\text{trace}([\|{\bf D}_1{\bf \Sigma}_1\|_2^2{\bf I} + {\bf \Sigma}_1^{-2}]^{-1})}\\
= & \|{\bf D}_1{\bf \Sigma}_1\|_2 \times \\
& \sqrt{\text{trace}({\bf \Sigma}_1[\|{\bf D}_1{\bf \Sigma}_1\|_2^2{\bf I} + {\bf \Sigma}_1^{-2}]^{-1}){\bf \Sigma}_1^T)}\\
\le & \dfrac{\sqrt{k}\|{\bf D}_1{\bf \Sigma}_1\|_2 \sigma_1} {\sqrt{{\sigma_1}^2 + \|{\bf D}_1{\bf \Sigma}_1\|_2^2}}.
\end{aligned} 
\end{equation}

By plugging this into \eqref{eqcomL1P1_second}, we obtain 
\begin{equation}
\begin{aligned}
\|({\bf I} - {\bf Q}_1 {\bf Q}_1^T){\bf A}_k\|_F  
\le & \dfrac{\sqrt{k}\|{\bf D}_1{\bf \Sigma}_1\|_2 \sigma_1}
{\sqrt{{\sigma_1}^2 + \|{\bf D}_1{\bf \Sigma}_1\|_2^2}}.
\end{aligned}
\label{eqIneqQ1bar}
\end{equation}

We are now prepared to present the result.
\begin{theorem} 
With the notation of Theorem \ref{Thr1}, and $\varrho=2, F$, the approximation error for Alg. \ref{Alg_first} must satisfy
\begin{equation}
\begin{aligned}
\|{\bf A} - \hat{\bf A}_\text{SOR}\|_\varrho \le 
& \|{{\bf A}_0}\|_\varrho +
\sqrt{\frac{\alpha^2{\|{{\bf \Omega}_2}\|_2^2} \|{{\bf \Omega}_1}^\dagger\|_2^2}{1 + \beta^2{\|{{\bf \Omega}_2}\|_2^2} \|{{\bf \Omega}_1}^\dagger\|_2^2}} \\ & + \sqrt{\frac{\eta^2{\|{{\bf \Omega}_2}\|_2^2} \|{{\bf \Omega}_1}^\dagger\|_2^2}{1 + \tau^2{\|{{\bf \Omega}_2}\|_2^2} \|{{\bf \Omega}_1}^\dagger\|_2^2}},
\end{aligned} 
\label{Thr2_1st}
\end{equation}
where $\alpha = \sqrt{k}\frac{\sigma_{\ell-p+1}^2}{\sigma_k}$, $\beta = \frac{\sigma_{\ell-p+1}^2}{\sigma_1 \sigma_k}$, $\eta=\sqrt{k}\sigma_{\ell-p+1}$ and $\tau=\frac{\sigma_{\ell-p+1}}{\sigma_1} $. 

When the power method is used, Alg. \ref{Alg3}, $\alpha = \sqrt{k}\frac{\sigma_{\ell-p+1}^2}{\sigma_k}\Big(\frac{\sigma_{\ell-p+1}}{\sigma_k} \Big)^{2q}$, $\beta = \frac{\sigma_{\ell-p+1}^2}{\sigma_1\sigma_k}\Big(\frac{\sigma_{\ell-p+1}}{\sigma_k} \Big)^{2q}$, $\eta=\dfrac{\sigma_k}{\sigma_{\ell-p+1}}\alpha$ and  $\tau=\dfrac{1}{\sigma_{\ell-p+1}}\beta$.
\label{Thr2}
\end{theorem} 

\textit{Proof.}
First, by plugging \eqref{eqcomL1P1_WSigma} and \eqref{eqIneqQ1bar} into 
\eqref{equThrFirst2} and \eqref{equThrFirst3}, for $\varrho=2, F$, we obtain 
\begin{equation}
\begin{aligned}
\|{\bf A} - \hat{\bf A}_\text{SOR}\|_\varrho \le &
\|{\bf A}_0\|_\varrho + \dfrac{\sqrt{k}\|{\bf W}_1{\bf \Sigma}_1\|_2 \sigma_1} {\sqrt{{\sigma_1}^2 + \|{\bf W}_1{\bf \Sigma}_1\|_2^2}} \\
& + \dfrac{\sqrt{k}\|{\bf D}_1{\bf \Sigma}_1\|_2 \sigma_1}
{\sqrt{{\sigma_1}^2 + \|{\bf D}_1{\bf \Sigma}_1\|_2^2}}.
\end{aligned}
\label{equThr2p_1}
\end{equation}

When the basic form of SOR-SVD algorithm is implemented, we can write 
${\bf W}_1{\bf \Sigma}_1$ as 
\begin{equation}
\begin{aligned}
{\bf W}_1{\bf \Sigma}_1 = {\bf \Sigma}_3^2{\bf \Omega}_2
{\bf \Omega}_1^\dagger {\bf \Sigma}_1^{-1}.
\end{aligned}
\end{equation}
where ${\bf W}_1$ is defined in \eqref{eqAtT1X}. Thus
\begin{equation}
\begin{aligned}
\|{\bf W}_1{\bf \Sigma}_1\|_2 \le \Big(\frac{\sigma_{\ell-p+1}^2}{\sigma_k}\Big)
\|{\bf \Omega}_2\|_2\|{\bf \Omega}_1^\dagger\|_2.
\end{aligned}
\label{eqWSigma}
\end{equation}

For ${\bf D}_1{\bf \Sigma}_1$, we can write 
\begin{equation}
\begin{aligned}
{\bf D}_1{\bf \Sigma}_1 = {\bf \Sigma}_3{\bf \Omega}_2{\bf \Omega}_1^\dagger.
\end{aligned}
\end{equation}
where ${\bf D}_1$ is defined in \eqref{eqD1_NoPM}. Thus
\begin{equation}
\begin{aligned}
\|{\bf D}_1{\bf \Sigma}_1\|_2 \le {\sigma_{\ell-p+1}} \|{\bf \Omega}_2\|_2\|{\bf \Omega}_1^\dagger\|_2.
\end{aligned}
\label{eqDSigma}
\end{equation}
 
Plugging \eqref{eqWSigma} and \eqref{eqDSigma} into \eqref{equThr2p_1}, 
and dividing both the numerator and denominator by $\sigma_1$, gives the result in \eqref{Thr2_1st}.

For the case when the power method is incorporated into the algorithm, by 
using \eqref{eqW1_PM} we can write ${\bf W}_1{\bf \Sigma}_1$ as
\begin{equation}
{\bf W}_1{\bf \Sigma}_1 = {\bf \Sigma}_3^{2q+2}{\bf \Omega}_2{\bf \Omega}_1^\dagger {\bf \Sigma}_1^{-(2q+1)}.
\end{equation} 

Consequently, we have
\begin{equation}
\begin{aligned}
\|{\bf W}_1{\bf \Sigma}_1\|_2 \le \frac{\sigma_{\ell-p+1}^2}{\sigma_k}\Big(\frac{\sigma_{\ell-p+1}}{\sigma_k}\Big)^{2q} \|{\bf \Omega}_2\|_2\|{\bf \Omega}_1^\dagger\|_2.
\end{aligned}
\label{eqWSigma2}
\end{equation}

To get ${\bf D}_1{\bf \Sigma}_1$, we first need to obtain ${\bf D}_1$, equation \eqref{eqD1_NoPM}, for the case when the power method is employed. To this end, the procedure starts with substituting ${\bf T}_1$ in equation \eqref{eqT1_2} with  
\begin{equation}
{\bf T}_1 = ({\bf AA}^T)^q{\bf A\Omega}.
\end{equation}

By writing the SVD of $\bf A$ \eqref{eqSVD}, we have
\begin{equation}
\begin{aligned} 
&{\bf T}_1 = {\bf U \Sigma}^{2q+1}{\bf V}^T \bf \Omega\\
& ={\bf U}\begin{bmatrix}
       {\bf \Sigma}_1^{2q+1} & {\bf 0} & {\bf 0} \\
       {\bf 0} & {\bf \Sigma}_2^{2q+1} & {\bf 0}\\
       {\bf 0} & {\bf 0} & {\bf \Sigma}_3^{2q+1}
  \end{bmatrix}
  \begin{bmatrix}
       {\bf V}_1^T{\bf \Omega} \\
       \\
       {\bf V}_2^T {\bf \Omega}
       \end{bmatrix} = {\bf Q}_1{\bf R}_1.
\end{aligned}
\end{equation}
and consequently $\bf X$ is defined as: 
\begin{equation}
{\bf X} = \bigg[
{\bf \Omega}_1^\dagger
   \begin{pmatrix}
        {\bf \Sigma}_1^{2q+1} & {\bf 0}\\
        {\bf 0} & {\bf \Sigma}_2^{2q+1}
   \end{pmatrix}^{-1}, \tilde{\bf X}\bigg].
\end{equation} 

Forming ${\bf T}_1{\bf X}$ yields:
\begin{equation}
{\bf D}_1 = {\bf \Sigma}_3^{2q+2}{\bf \Omega}_2{\bf \Omega}_1^\dagger {\bf \Sigma}_1^{-(2q+1)}.
\label{eqD1_PM}
\end{equation}
and we have
\begin{equation}
\begin{aligned}
{\bf D}_1{\bf \Sigma}_1 = {\bf \Sigma}_3^{2q+1}{\bf \Omega}_2{\bf \Omega}_1^\dagger {\bf \Sigma}_1^{-2q}.
\end{aligned}
\end{equation}

Accordingly, we have
\begin{equation}
\begin{aligned}
\|{\bf D}_1{\bf \Sigma}_1\|_2 \le {\sigma_{\ell-p+1}}\Big(\frac{\sigma_{\ell-p+1}}{\sigma_k}\Big)^{2q}
\|{\bf \Omega}_2\|_2\|{\bf \Omega}_1^\dagger\|_2.
\end{aligned}
\label{eqDSigma2}
\end{equation}

Substituting \eqref{eqWSigma2} and \eqref{eqDSigma2} into \eqref{equThr2p_1}, 
and dividing both the numerator and denominator by $\sigma_1$, gives the result. \QEDB

Theorem \ref{Thr1} shows that the accuracy of singular values depends strongly on the ratio $\frac{\sigma_{l-p+1}}{\sigma_j}$ for $j = 1,..., k$, whereas by Theorem
\ref{Thr2}, the accuracy of the low-rank approximation depends on $\frac{\sigma_{l-p+1}}{\sigma_k}$. The power method decreases the extra factors in the error bounds by driving down the aforesaid ratios exponentially fast. Thus, by increasing the number of iterations $q$, we can make the extra factors as close to zero as we wish. However, this increases the cost of the algorithm.

\subsection{Average-Case Error Bounds}
\label{secProEB}

In this section, we provide an average-case error analysis for the SOR-SVD 
algorithm, which, in contrast to the argument in Section \ref{secDetEB}, 
depends on distributional assumptions on the random matrix $\bf \Omega$. To 
be precise, $\bf \Omega$ has a standard Gaussian distribution which is invariant under all orthogonal transformations. This means if, for instance, $\bf V$ has orthonormal columns, then the product ${\bf V}^T{\bf \Omega}$ has the same standard Gaussian distribution. This allows us to take advantage of the vast 
literature on properties of Gaussian matrices.

We begin by stating a few propositions that are used later on. 

\begin{proposition}(Gu \cite{Gu2015}) Let ${\bf G}\in \mathbb R^{m \times n}$ be a 
Gaussian matrix. For any $\alpha>0$ we have
\begin{equation}
\mathbb{E}\Bigg(\frac{1}{\sqrt{1 + \alpha^2\|{\bf G}\|_2^2}}\Bigg) \ge  
\frac{1}{\sqrt{1 + \alpha^2\nu^2}},
\label{eqGu_Gsv}
\end{equation}
where $\nu = \sqrt{m}+\sqrt{n}+7$.
\label{Gu_Gsv}
\end{proposition}

\begin{proposition} With the notation of Proposition \ref{Gu_Gsv}, and furthermore, for any $\beta>0$, we have
\begin{equation}
\mathbb{E}\Bigg(\sqrt{\frac{\alpha^2\|{\bf G}\|_2^2}{1 + \beta^2\|{\bf G}\|_2^2}}\Bigg) \le \sqrt{\frac{\alpha^2\nu^2}{1 + \beta^2\nu^2}}.
\label{eqProp_G} 
\end{equation}
\label{Prop_G} 
where $\nu$ is defined in \eqref{eqGu_Gsv}. 
\end{proposition}

\textit{Proof.} The proof is given in Appendix \ref{secAProp_G}.

\begin{proposition}(Gu \cite{Gu2015}) Let ${\bf G}\in \mathbb R^{\ell-p \times \ell}$ be a Gaussian matrix. Then $\text{rank}({\bf G}) = \ell-p$ with probability 1. For $p \ge 2$ and any $\alpha>0$ 
\begin{equation}
\mathbb{E}\Bigg(\frac{1}{\sqrt{1 + \alpha^2\|{\bf G}^\dagger\|_2^2}}\Bigg) \ge  
\frac{1}{\sqrt{1 + \alpha^2\nu^2}},
\label{eqGu_PIsv}
\end{equation}
where $\nu = \dfrac{4\text{e}\sqrt{\ell}}{p+1}$. 
\label{Gu_PIsv}
\end{proposition} 
 
\begin{proposition} With the notation of Proposition \ref{Gu_PIsv}, and 
furthermore, for any $\beta>0$, we have
\begin{equation}
\mathbb{E}\Bigg(\sqrt{\frac{\alpha^2\|{\bf G}^\dagger\|_2^2}{1 + \beta^2\|{\bf G}^\dagger\|_2^2}}\Bigg) \le \sqrt{\frac{\alpha^2\nu^2}{1 + \beta^2\nu^2}}.
\end{equation}
\label{Prop_PI}
where $\nu$ is defined in \eqref{eqGu_PIsv}.
\end{proposition}

\textit{Proof.} The proof is given in Appendix \ref{secAProp_PI}.

We now present the result.

\begin{theorem}
With the notation of Theorem \ref{Thr1}, and $\gamma_j = \frac{\sigma_{\ell-p+1}}{\sigma_j}$, for $j = 1,..., k,$, for Alg. \ref{Alg_first}, we have  
\begin{equation}
\mathbb{E}(\sigma_j(\hat{\bf A}_\text{SOR})) \ge \frac{\sigma_j}
{\sqrt{1 + \nu^2\gamma_j^4}},
\end{equation}
and when the power method is used, Alg. \ref{Alg3}, we have 
\begin{equation}
\begin{aligned}
\mathbb{E}(\sigma_j(\hat{\bf A}_\text{SOR})) \ge \frac{\sigma_j}
{\sqrt{1 + \nu^2\gamma_j^{4q+4}}},
\end{aligned}
\label{eqThr3_2nd}
\end{equation}

where $\nu = \nu_1\nu_2$, $\nu_1 = \sqrt{n-\ell+p}+\sqrt{\ell}+7$, and $\nu_2 = \frac{4\text{e}\sqrt{\ell}}{p+1}$.
\label{Thr3}
\end{theorem}

\textit{Proof.}
We only prove Theorem \ref{Thr3} for $j = k$, as is the case for Theorem \ref{Thr1}. Theorem \ref{Thr3} can be proved for other values of $1 \le j<k$ by referring to Theorem \ref{Thr3} for a rank $j$ truncated SVD. Since ${\bf \Omega}_1$ and ${\bf \Omega}_2$ are independent from each other, to bound expectations, we  in turn take expectations over ${\bf \Omega}_2$ and ${\bf \Omega}_1$:\looseness-1
\begin{equation}
\begin{aligned}
& \mathbb{E}(\sigma_k(\hat{\bf A}_\text{SOR})) = 
\mathbb{E}_{{\bf \Omega}_1}  \big(\mathbb{E}_{{\bf \Omega}_2} \big[\sigma_k(\hat{\bf A}_\text{SOR})\big]\big)\\ 
& = \mathbb{E}_{{\bf \Omega}_1}  \left(\mathbb{E}_{{\bf \Omega}_2}\left[\frac{\sigma_j}{\sqrt{1 + {\|{{\bf \Omega}_2}\|_2^2} \|{{\bf \Omega}_1}^\dagger\|_2^2 \Big(\frac{\sigma_{\ell-p+1}}{\sigma_j}\Big)^4}}\right]\right) \\
& \ge \mathbb{E}_{{\bf \Omega}_1} \left(\frac{\sigma_k}{\sqrt{1 + \nu_1^2\|{{\bf \Omega}_1}^\dagger\|_2^2 \gamma_j^4}}\right)\\ 
& \ge \frac{\sigma_k}{\sqrt{1 + \nu^2\gamma_k^4}}.
\end{aligned}
\label{eqThr3_proof}
\end{equation}
                
The second line follows from Theorem \ref{Thr1}, equation \eqref{eqThr1_sv1}, the third line from Proposition \ref{Gu_Gsv}, and the last line from Proposition \ref{Gu_PIsv}. The result in 
\eqref{eqThr3_2nd}, likewise, follows by substituting \eqref{eqThr1_sv2} into the second line of equation \eqref{eqThr3_proof}. \QEDB

We now present a theorem that establishes average error bounds on the low-rank
approximation.

\begin{theorem}
With the notation of Theorem \ref{Thr1}, and $\varrho = 2,F$, for Alg. \ref{Alg_first}, we have  
\begin{equation}
\begin{aligned}
\mathbb{E} \|{{\bf A} - \hat{\bf A}_\text{SOR}}\|_\varrho \le \mspace{5mu} \|{{\bf A}_0}\|_\varrho + 
(1 + \gamma_k)
\sqrt{k}\nu\sigma_{\ell-p+1},
\end{aligned} 
\end{equation}
and when the power method is used, Alg. \ref{Alg3}, we have
\begin{equation}
\begin{aligned}
\mathbb{E} \|{{\bf A} - \hat{\bf A}_\text{SOR}}\|_\varrho \le \mspace{5mu} \|{{\bf A}_0}\|_\varrho + (1 + \gamma_k) \sqrt{k}\nu\sigma_{\ell-p+1}\gamma_k^{2q},
\end{aligned} 
\end{equation}
where $\gamma_k$ and $\nu$ are defined in Theorem \ref{Thr3}.
\label{Thr4}
\end{theorem}

\textit{Proof.} Similar to the proof of Theorem \ref{Thr3}, we first take expectations over ${\bf \Omega}_2$ and next over ${\bf \Omega}_1$. By invoking
Theorem \ref{Thr2}, Propositions \ref{Prop_G}, and Proposition \ref{Prop_PI}, 
we have 

\begin{equation}
\begin{aligned}
& \mathbb{E}\|{\bf A} - \hat{\bf A}_\text{SOR}\|_\varrho = 
\mathbb{E}_{{\bf \Omega}_1} \mathbb{E}_{{\bf \Omega}_2}\|{\bf A} - \hat{\bf A}_\text{SOR}\|_\varrho\\
& = \|{{\bf A}_0}\|_\varrho + \mathbb{E}_{{\bf \Omega}_1} \mathbb{E}_{{\bf \Omega}_2} \Bigg(
\sqrt{\frac{\alpha^2{\|{{\bf \Omega}_2}\|_2^2} \|{{\bf \Omega}_1}^\dagger\|_2^2}{1 + \beta^2{\|{{\bf \Omega}_2}\|_2^2} \|{{\bf \Omega}_1}^\dagger\|_2^2}} \\ & + \sqrt{\frac{\eta^2{\|{{\bf \Omega}_2}\|_2^2} \|{{\bf \Omega}_1}^\dagger\|_2^2}{1 + \tau^2{\|{{\bf \Omega}_2}\|_2^2} \|{{\bf \Omega}_1}^\dagger\|_2^2}} \Bigg)\\ 
& \le \|{{\bf A}_0}\|_\varrho + \mathbb{E}_{{\bf \Omega}_1}
\Bigg(\sqrt{\frac{\alpha^2{\nu_1^2} \|{{\bf \Omega}_1}^\dagger\|_2^2}{1 + \beta^2{\nu_1^2} \|{{\bf \Omega}_1}^\dagger\|_2^2}} \\ & + 
\sqrt{\frac{\eta^2{\nu_1^2} \|{{\bf \Omega}_1}^\dagger\|_2^2}{1 + \tau^2{\nu_1^2} \|{{\bf \Omega}_1}^\dagger\|_2^2}}\Bigg)\\
& \le \|{{\bf A}_0}\|_\varrho + \sqrt{\frac{\alpha^2{\nu_1^2}\nu_2^2}{1 + \beta^2{\nu_1^2}\nu_2^2}}  + 
 \sqrt{\frac{\eta^2{\nu_1^2} \nu_2^2}{1 + \tau^2{\nu_1^2}\nu_2^2}}\\
& \le \|{{\bf A}_0}\|_\varrho + (\alpha + \eta)\nu.
\end{aligned}
\end{equation}

Plugging values of $\alpha$, $\eta$, defined in Theorem \ref{Thr2}, and $\nu$ into the last inequality gives 
the results. \QEDB

\subsection{Computational Complexity}
\label{secComComplex}
The cost of any algorithm involves both arithmetic, i.e., the number of floating-point operations, and communication, i.e., synchronization and data movement either through levels of a memory hierarchy or between parallel processing units \cite{DemmGHL12}. On multicore and accelerator-based computers, for a data matrix which is too large to fit in fast memory, the communication cost becomes substantially more expensive compared to the arithmetic \cite{DemmGHL12, Dongarra17}. Thus, developing new algorithms, or redesigning existing algorithms, 
to solve a problem in hand with minimum communication cost is highly desirable.\looseness-1

The advantage of randomized algorithms over their classical counterparts lies in the fact that \textit{i)} they operate on a compressed version of the data matrix rather than a matrix itself, and \textit{ii)} they can be organized to take advantage of modern architectures, performing a decomposition with minimum communication cost. 

To decompose the matrix ${\bf A}$, the simple version of the SOR-SVD algorithm, requires only three passes over data, for a matrix stored-out-of-core, with the following operation count:\looseness-1
\begin{equation}
C_\text{SOR-SVD} \sim 3\ell C_\text{mult}+2\ell (m+n)(\ell+k)+2l^2(m+\ell),
\label{equC1}
\end{equation}
where $C_\text{mult}$ is the cost of a matrix-vector multiplication with 
$\bf A$ or ${\bf A}^T$, and the cost of a QR factorization or an SVD, for 
instance, for $\bf A$, is $2mn^2$ flops \cite{GolubVanLoan96}. The flop 
count for the two-pass SOR-SVD, where the matrix $\bf M$ is approximated 
by ${\bf M}_\text{approx}$, satisfies
\begin{equation}
C_\text{SOR-SVD} \sim 2\ell C_\text{mult}+2\ell (m+n)(2\ell + k)+5\ell ^3.
\label{equC1_2pass}
\end{equation}

When the power method is employed, Alg. \ref{Alg3}, the algorithm require either $2q+3$ passes through data with the following 
operation count:
\begin{equation}
C_\text{SOR-SVD+PM} \sim (2q+3)\ell C_\text{mult}+2\ell (m+n)(\ell+k)+2\ell^2(m+\ell),
\label{equC2}
\end{equation}
or $2q+2$ passes over data with the following flop count:
\begin{equation}
C_\text{SOR-SVD+PM} \sim (2q+2)\ell C_\text{mult}+2\ell (m+n)(2\ell +k)+5\ell ^3.
\label{equC2_2pass}
\end{equation}

Except for matrix-matrix multiplications, which are easily parallelizable, the SOR-SVD algorithm performs two QR decompositions on matrices of size $m\times \ell$ and $n\times \ell$, and one SVD on an $\ell \times \ell$ matrix, whereas the R-SVD performs one QR decomposition on an $m\times \ell$ matrix and one SVD on a $n\times \ell$ matrix. While recently developed Communication-Avoiding QR (CAQR) algorithms \cite{DemmGHL12} are optimal in terms of communication costs, standard techniques to compute an SVD are challenging for parallelization \cite{Demmel97, MarQH17}. Having known that for a general large $\bf A$, $k \ll \text{min}\{m,n\}$, and moreover, today's computing devices have hardware switches being controlled in software \cite{Dongarra17}, the SVD of the $\ell \times \ell$ matrix can be computed 
either \textit{in-core} on a sequential processor, or with minimum communication cost on distributed processors. This significantly reduces the computational time of SOR-SVD, a superiority of the algorithm over R-SVD.

As explained, degradation of the quality of approximation by TSR-SVD, Alg. \ref{Alg2}, as a result of one pass over the data, makes it unsuitable for use in practice. Moreover, as will be shown in the next section, the SOR-SVD algorithm shows even better performance compared to the two-pass (or $2q+2$ passes when the power method is used) TSR-SVD algorithm (this is consistent with the theoretical bounds obtained).

\section{Numerical Experiments}
\label{secNumExp}

In this section, we evaluate the performance of the SOR-SVD algorithm through numerical tests. Our goal is to experimentally investigate the behavior of the SOR-SVD algorithm in different scenarios.

In Section \ref{SynMat}, we consider two classes of synthetic matrices to experimentally verify that the SOR-SVD algorithm provides highly accurate singular values and low-rank approximations. We compare the performance of the SOR-SVD with those of the SVD, R-SVD, and TSR-SVD.

In Section \ref{subEmpEva}, we experimentally investigate the tightness of the low-rank approximation error bounds for the spectral and Frobenius norm provided 
in Theorem \ref{Thr2}.

In Section \ref{Exrpca}, we apply the SOR-SVD to recovering a low-rank plus a sparse matrix through experiments on \textit{i}) synthetically generated data with various dimensions, numerical rank and gross errors, and \textit{ii}) real-time data in applications to background subtraction in surveillance video, and shadow and specularity removal from face images. 

\subsection{Synthetic Matrices}
\label{SynMat}
We first describe two types of data matrices that we use in our tests to assess the behavior of the SOR-SVD algorithm. For the sake of simplicity, we focus on square matrices. 
\subsubsection{A Noisy Low-rank Matrix}
\label{NoisyLRM}
A rank-$k$ matrix ${\bf A}$ of order $1000$ is generated as ${\bf A} ={\bf U\Sigma V}^T +0.1\sigma_k{\bf E}$, introduced by Stewart \cite{StewartQLP}, where ${\bf U}$ and ${\bf V}$ are random orthonormal matrices, ${\bf \Sigma}$ is diagonal containing the singular values $\sigma_i$s that decrease geometrically from $1$ to $10^{-9}$, and $\sigma_{j}=0$, for $j=k+1, ..., 1000$. The matrix ${\bf E}$ is a normalized Gaussian matrix. Here, we set $k=20$.

\subsubsection{A Matrix with Rapidly Decaying Singular Values}
A matrix ${\bf A}$ of order $n$ is generated as ${\bf A} ={\bf U\Sigma V}^T$, where 
${\bf U}$ and ${\bf V}$ are random orthonormal matrices, ${\bf \Sigma}= \text{diag}(1^{-1}, 2^{-1}, ..., n^{-1})$. We set $n=10^3$, and the rank $k=10$.  

To study the behavior of the SOR-SVD, and to provide a fair comparison, we 
factor the matrix $\bf A$ using the SVD, R-SVD, TSR-SVD, and SOR-SVD algorithms.
For the randomized algorithms, we set sample size parameter $\ell=38$ for the 
first matrix, and $\ell=18$ for the second one, chosen randomly. All randomized algorithms require the same number of passes through $\bf A$, either two or 
$2q+2$ when the power iteration is used. We then compare the singular values 
and low-rank approximations computed by the algorithms mentioned.

\begin{figure}[t]
\begin{center}
\input{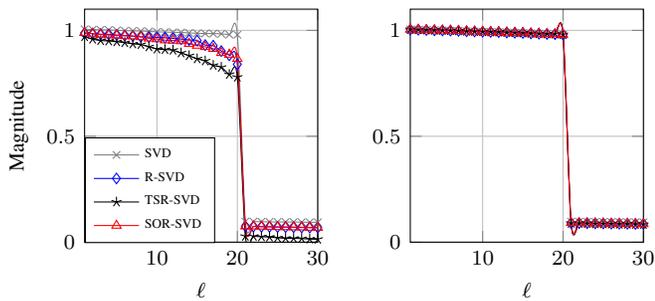}
\captionsetup{justification=centering,font=scriptsize}
\caption{Comparison of singular values for the noisy low-rank matrix. No power method, $(q=0)$ (left), and $q=2$ (right).} 
\label{fig:SinValu_1}       
\end{center}
\end{figure}

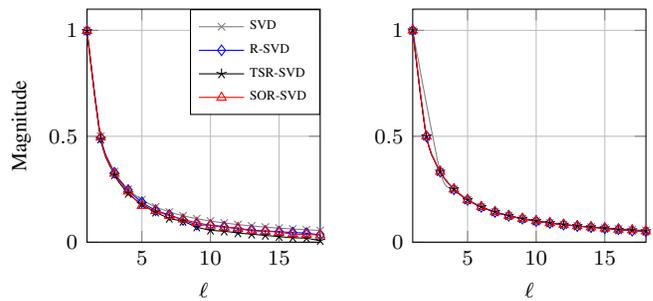
\begin{figure}[t]
\begin{center}
%
%
%
\usetikzlibrary{positioning,calc}

\definecolor{mycolor1}{rgb}{0.00000,1.00000,1.00000}%
\definecolor{mycolor2}{rgb}{1.00000,0.00000,1.00000}%

\pgfplotsset{every axis label/.append style={font=\footnotesize},
every tick label/.append style={font=\footnotesize}
}

\begin{tikzpicture}[font=\footnotesize] 

\begin{axis}[%
name=ber,
width  = 0.35\columnwidth,
height = 0.35\columnwidth,
scale only axis,
xmin  = 1,
xmax  = 18,
xlabel= {$\ell$},
xmajorgrids,
ymin=0,
ymax=1.1,
ylabel={Magnitude},
ymajorgrids,
legend entries={SVD, R-SVD, TSR-SVD,SOR-SVD},
legend style={at={(1,1)},anchor=north east,draw=black,fill=white,legend cell align=left,font=\tiny}
]


\addplot+[smooth,color=gray,solid, every mark/.append style={solid}, mark=x]
table[row sep=crcr]{
1	0.999999999999999\\
2	0.500000000000001\\
3	0.333333333333335\\
4	0.250000000000000\\
5	0.199999999999999\\
6	0.166666666666667\\
7	0.142857142857143\\
8	0.125000000000001\\
9	0.111111111111112\\
10	0.100000000000000\\
11	0.0909090909090911\\
12	0.0833333333333335\\
13	0.0769230769230771\\
14	0.0714285714285717\\
15	0.0666666666666666\\
16	0.0624999999999999\\
17	0.0588235294117650\\
18	0.0555555555555557\\
19	0.0526315789473687\\
20	0.0500000000000002\\
}; 

\addplot+[smooth,color=blue,solid, every mark/.append style={solid}, mark = diamond]
  table[row sep=crcr]{
1	0.998386592368603\\
2	0.493631887628756\\
3	0.328187647028117\\
4	0.246327536487902\\
5	0.188397710470007\\
6	0.151425407365790\\
7	0.128918490343693\\
8	0.105232739421815\\
9	0.0834077545503760\\
10	0.0802898667435362\\
11	0.0732943056403757\\
12	0.0628534504848684\\
13	0.0572358459749275\\
14	0.0526821164018446\\
15	0.0500967296634081\\
16	0.0452975775457965\\
17	0.0428648175048246\\
18	0.0342300631990984\\
19	0.0293582104162083\\
20	0.0253218488140999\\
};

\addplot+[smooth,color=black ,solid, every mark/.append style={solid}, mark=star]
  table[row sep=crcr]{
1	0.996704091641284\\
2	0.487256530244323\\
3	0.317804766101812\\
4	0.230505543914299\\
5	0.180345162576855\\
6	0.141918114206955\\
7	0.112775553606834\\
8	0.0983644530850231\\
9	0.0718912307823512\\
10	0.0564979406503854\\
11	0.0519633431153076\\
12	0.0436938823629285\\
13	0.0378494241784325\\
14	0.0322259726090601\\
15	0.0263816039719368\\
16	0.0207752301594806\\
17	0.0187769350026855\\
18	0.00753322195971780\\
19	0.00538223204547392\\
20	0.00132792153602228\\
};

\addplot+[smooth,color=red, solid, every mark/.append style={solid}, mark=triangle]
  table[row sep=crcr]{
1	0.997854515104833\\
2	0.495078787606555\\
3	0.324241214166147\\
4	0.239935245609783\\
5	0.172116815493796\\
6	0.148340399210366\\
7	0.124968690353354\\
8	0.108206312629543\\
9	0.0916160889699227\\
10	0.0806277894736536\\
11	0.0698650675062045\\
12	0.0677061895638335\\
13	0.0538203513838833\\
14	0.0518525924730861\\
15	0.0471503471619046\\
16	0.0382623462412133\\
17	0.0372386559440389\\
18	0.0331971213536172\\
19	0.0304345473798118\\
20	0.0275584620139238\\
}; 

\end{axis}

\begin{axis}[%
name=SumRate,
at={($(ber.east)+(35,0em)$)},
		anchor= west,
width  = 0.35\columnwidth,
height = 0.35\columnwidth,
scale only axis,
xmin  = 1,
xmax  = 18,
xlabel= {$\ell$},
xmajorgrids,
ymin=0,
ymax=1.1,
ylabel={},
ymajorgrids,
]


\addplot+[smooth,color=gray,solid, every mark/.append style={solid}, mark = x]
table[row sep=crcr]{
1	1
2	0.500000000000002\\
3	0.333333333333333\\
4	0.249999999999999\\
5	0.200000000000000\\
6	0.166666666666667\\
7	0.142857142857143\\
8	0.125000000000000\\
9	0.111111111111111\\
10	0.100000000000001\\
11	0.0909090909090913\\
12	0.0833333333333332\\
13	0.0769230769230772\\
14	0.0714285714285717\\
15	0.0666666666666668\\
16	0.0624999999999999\\
17	0.0588235294117649\\
18	0.0555555555555552\\
19	0.0526315789473680\\
20	0.0499999999999995\\
};

\addplot+[smooth,color = blue, solid, every mark/.append style={solid}, mark=diamond]
  table[row sep=crcr]{
1	0.999999999999976\\
2	0.499999999987337\\
3	0.333333333183799\\
4	0.249999992117739\\
5	0.199999972224774\\
6	0.166666566295947\\
7	0.142853973670104\\
8	0.124997626451210\\
9	0.111102963366083\\
10	0.0999870972783954\\
11	0.0906928968217901\\
12	0.0832039711777330\\
13	0.0762680408300673\\
14	0.0712440695753269\\
15	0.0660412729262119\\
16	0.0614490736095275\\
17	0.0546802863828324\\
18	0.0513927575868093\\
19	0.0495064684793235\\
20	0.0491171947012105\\
};

\addplot+[smooth,color=black,solid, every mark/.append style={solid}, mark=star]
  table[row sep=crcr]{
1	0.999999999990317\\
2	0.499999998190555\\
3	0.333333317187662\\
4	0.249999855209577\\
5	0.199999474759919\\
6	0.166652140218610\\
7	0.142835214753726\\
8	0.124983642633173\\
9	0.111036329730273\\
10	0.0999064580039391\\
11	0.0908515765326129\\
12	0.0829959455297590\\
13	0.0763872460252986\\
14	0.0709020194999323\\
15	0.0660406020122904\\
16	0.0585801718910436\\
17	0.0567369998495021\\
18	0.0521272733865546\\
19	0.0481387324993503\\
20	0.0429029082556785\\
};

\addplot+[smooth,color=red, solid, every mark/.append style={solid}, mark=triangle]
  table[row sep=crcr]{
1	0.999999999999982\\
2	0.499999999940439\\
3	0.333333332918173\\
4	0.249999982123241\\
5	0.199999900467244\\
6	0.166665844346940\\
7	0.142853997903128\\
8	0.124998083611501\\
9	0.111097765555481\\
10	0.0999777008555389\\
11	0.0908502457557932\\
12	0.0831287460457059\\
13	0.0767653675724008\\
14	0.0712172627717329\\
15	0.0657264826149644\\
16	0.0603662125731166\\
17	0.0575850436571018\\
18	0.0523196708315239\\
19	0.0494863173321897\\
20	0.0450332090632578\\
};

\end{axis}

\end{tikzpicture}%
\captionsetup{justification=centering,font=scriptsize}
\caption{Comparison of singular values for the matrix with polynomially decaying singular values. No power method, $q=0$, (left) and $q=2$ (right).} 
\label{fig:SinValu_2}       
\end{center}
\end{figure}

Figures \ref{fig:SinValu_1} and \ref{fig:SinValu_2} display the singular values computed by each algorithm for the two classes of matrices. The SOR-SVD uses a truncated SVD on the reduced-size matrix $\bf M$, however, we show the results 
for the full SVD, i.e., a full SVD of the $\ell \times \ell$ matrix, for the 
purpose of comparison.  

Judging from the figures, when the power method is not used, 
$q = 0$, which for the SOR-SVD corresponds to the basic form of the algorithm, 
the SOR-SVD approximates for some singular values are good estimates of the 
true ones, while for the others are relatively poor. When we set $q=2$, we 
see a substantial improvement in the quality of the approximations.

Nevertheless, we observe that when $q=0$, the SOR-SVD approximates of 
both leading and trailing singular values exceed those of the TSR-SVD, 
while being similar to those of the R-SVD. In the case when $q=2$, as 
shown, all methods perform similarly.

To compare the quality of low-rank approximations, we first compute a rank-$k$ approximation ${\hat{\bf A}}_\text{out}$ to ${\bf A}$ by each algorithm. We 
do so by varying the sample size parameter $\ell$, while the rank is fixed. We then calculate the Frobenius-norm error:
\begin{equation}
e_k = \|{\bf A} - \hat{\bf A}_{\text{out}}\|_F.
\end{equation}

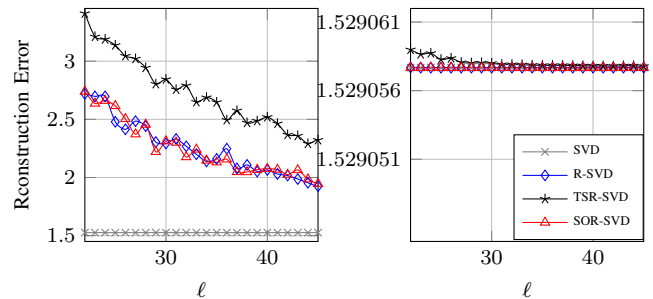
\begin{figure}[t]
\begin{center}
%
%
%
\usetikzlibrary{positioning,calc}

\definecolor{mycolor1}{rgb}{0.00000,1.00000,1.00000}%
\definecolor{mycolor2}{rgb}{1.00000,0.00000,1.00000}%

\pgfplotsset{every axis label/.append style={font=\footnotesize},
every tick label/.append style={font=\footnotesize}
}

\begin{tikzpicture}[font=\footnotesize] 

\begin{axis}[%
name=ber,
width  = 0.35\columnwidth,
height = 0.35\columnwidth,
scale only axis,
xmin  = 22,
xmax  = 45,
xlabel= {$\ell$},
xmajorgrids,
ymin = 1.45,
ymax = 3.4533,
ylabel={Rconstruction Error},
ymajorgrids,
]


\addplot+[color=gray, mark=x]
table[row sep=crcr]{
22	1.52584568568731\\
23	1.52584568568731\\
24	1.52584568568731\\
25	1.52584568568731\\
26	1.52584568568731\\
27	1.52584568568731\\
28	1.52584568568731\\
29	1.52584568568731\\
30	1.52584568568731\\
31	1.52584568568731\\
32	1.52584568568731\\
33	1.52584568568731\\
34	1.52584568568731\\
35	1.52584568568731\\
36	1.52584568568731\\
37	1.52584568568731\\
38	1.52584568568731\\
39	1.52584568568731\\
40	1.52584568568731\\
41	1.52584568568731\\
42	1.52584568568731\\
43	1.52584568568731\\
44	1.52584568568731\\
45	1.52584568568731\\
46	1.52584568568731\\
47	1.52584568568731\\
48	1.52584568568731\\
49	1.52584568568731\\
50	1.52584568568731\\
}; 

\addplot+[smooth,color = blue,solid, every mark/.append style={solid}, mark=diamond]
  table[row sep=crcr]{
22	2.71924489918236\\
23	2.69712569855466\\
24	2.69706073781344\\
25	2.47882665318856\\
26	2.41395940646083\\
27	2.48446936162652\\
28	2.44119790074854\\
29	2.30461338560989\\
30	2.29249773374743\\
31	2.33152202015201\\
32	2.26808058516412\\
33	2.19950389007215\\
34	2.13821405866209\\
35	2.15987126174506\\
36	2.24731942680021\\
37	2.07770551892139\\
38	2.10850842310262\\
39	2.05105983866138\\
40	2.06462831528805\\
41	2.02994376400097\\
42	2.01647979715504\\
43	1.98638940204807\\
44	1.95482265495728\\
45	1.92484619659035\\
46	1.95320718305662\\
47	1.98472309616672\\
48	1.95945926956667\\
49	1.89913545289191\\
50	1.91119920747538\\
};

\addplot+[smooth,color=black,solid, every mark/.append style={solid}, mark=star]
  table[row sep=crcr]{
22	3.41042572583709\\
23	3.21047812513411\\
24	3.18786438620061\\
25	3.13673122165209\\
26	3.04462830738495\\
27	3.02037571373979\\
28	2.94311302711182\\
29	2.80344914731445\\
30	2.84365545721677\\
31	2.75602219276060\\
32	2.79020244806588\\
33	2.64880263251788\\
34	2.68783100271295\\
35	2.64650496013337\\
36	2.49386080495862\\
37	2.57466110806798\\
38	2.47014397774633\\
39	2.48581924998133\\
40	2.51856798744873\\
41	2.46292191681929\\
42	2.36662601700642\\
43	2.35768926060065\\
44	2.29026036501884\\
45	2.32169612055039\\
46	2.25864235307444\\
47	2.27268221462892\\
48	2.24575581944794\\
49	2.22129783653767\\
50	2.19460750036913\\
};

\addplot+[smooth,color=red,solid, every mark/.append style={solid}, mark=triangle]
  table[row sep=crcr]{
22	2.73882154295403\\
23	2.63632750502425\\
24	2.65991478542417\\
25	2.61560431192238\\
26	2.50367348483734\\
27	2.37226405707853\\
28	2.45541519035768\\
29	2.22112524928816\\
30	2.31233889676596\\
31	2.30016665672298\\
32	2.17568837206831\\
33	2.24123234300103\\
34	2.14679409225607\\
35	2.13135498965511\\
36	2.15641481809047\\
37	2.04778892521146\\
38	2.04691329941554\\
39	2.06746051191830\\
40	2.07542734711529\\
41	2.06801555445710\\
42	2.02495855804230\\
43	2.06430520074508\\
44	1.98535146314484\\
45	1.94673932503931\\
46	1.93995799690897\\
47	1.97595514543229\\
48	1.93890943345795\\
49	1.90096233784372\\
50	1.91215710250844\\
}; 

\end{axis}

\begin{axis}[%
name=SumRate,
at={($(ber.east)+(35,0em)$)},
		anchor= west,
width  = 0.35\columnwidth,
height = 0.35\columnwidth,
scale only axis,
xmin  = 22,
xmax  = 45,
xlabel= {$\ell$},
xmajorgrids,
ymin = 1.529045,
ymax = 1.529062,
ylabel={},
ymajorgrids,
ytick       ={1.529051, 1.529056, 1.529061},
yticklabels ={$1.529051$, $1.529056$ , $1.529061$},
legend entries={SVD, R-SVD, TSR-SVD,SOR-SVD},
legend style={at={(1,0.46)},anchor=north east,draw=black,fill=white,legend cell align=left,font=\tiny}
]

\addplot+[color=gray, mark=x]
table[row sep=crcr]{
22	1.52905767137784\\
23	1.52905767137784\\
24	1.52905767137784\\
25	1.52905767137784\\
26	1.52905767137784\\
27	1.52905767137784\\
28	1.52905767137784\\
29	1.52905767137784\\
30	1.52905767137784\\
31	1.52905767137784\\
32	1.52905767137784\\
33	1.52905767137784\\
34	1.52905767137784\\
35	1.52905767137784\\
36	1.52905767137784\\
37	1.52905767137784\\
38	1.52905767137784\\
39	1.52905767137784\\
40	1.52905767137784\\
41	1.52905767137784\\
42	1.52905767137784\\
43	1.52905767137784\\
44	1.52905767137784\\
45	1.52905767137784\\
46	1.52905767137784\\
47	1.52905767137784\\
48	1.52905767137784\\
49	1.52905767137784\\
50	1.52905767137784\\
};

\addplot+[smooth,color=blue,solid, every mark/.append style={solid}, mark=diamond]
  table[row sep=crcr]{
22	1.52905768813294\\
23	1.52905767803950\\
24	1.52905767769031\\
25	1.52905767686585\\
26	1.52905767596068\\
27	1.52905767464448\\
28	1.52905767449021\\
29	1.52905767354400\\
30	1.52905767437982\\
31	1.52905767361049\\
32	1.52905767362053\\
33	1.52905767305234\\
34	1.52905767275183\\
35	1.52905767290334\\
36	1.52905767252720\\
37	1.52905767245182\\
38	1.52905767243005\\
39	1.52905767229814\\
40	1.52905767237341\\
41	1.52905767223886\\
42	1.52905767218162\\
43	1.52905767215370\\
44	1.52905767218482\\
45	1.52905767210018\\
46	1.52905767198439\\
47	1.52905767212428\\
48	1.52905767209590\\
49	1.52905767205626\\
50	1.52905767196420\\
};

\addplot+[smooth,color=black, mark=star]
  table[row sep=crcr]{
22	1.52905898478907\\
23	1.52905865330401\\
24	1.52905872380971\\
25	1.52905827694757\\
26	1.52905836380777\\
27	1.52905805279444\\
28	1.52905800996423\\
29	1.52905802417820\\
30	1.52905799413027\\
31	1.52905790120027\\
32	1.52905790342030\\
33	1.52905789632975\\
34	1.52905786857281\\
35	1.52905784627828\\
36	1.52905784434682\\
37	1.52905783994738\\
38	1.52905782534475\\
39	1.52905782066806\\
40	1.52905780438682\\
41	1.52905779676501\\
42	1.52905779248571\\
43	1.52905782088246\\
44	1.52905781497083\\
45	1.52905777509052\\
46	1.52905778040116\\
47	1.52905777007522\\
48	1.52905774716159\\
49	1.52905776655842\\
50	1.52905777411950\\
};

\addplot+[smooth,color=red,solid, every mark/.append style={solid}, mark=triangle]
  table[row sep=crcr]{
22	1.52905767715806\\
23	1.52905767626272\\
24	1.52905768109457\\
25	1.52905767623574\\
26	1.52905767514640\\
27	1.52905767427735\\
28	1.52905767365289\\
29	1.52905767304147\\
30	1.52905767312373\\
31	1.52905767289541\\
32	1.52905767364936\\
33	1.52905767319115\\
34	1.52905767299195\\
35	1.52905767260769\\
36	1.52905767258164\\
37	1.52905767251926\\
38	1.52905767234861\\
39	1.52905767251171\\
40	1.52905767232023\\
41	1.52905767226893\\
42	1.52905767213574\\
43	1.52905767223514\\
44	1.52905767214279\\
45	1.52905767214744\\
46	1.52905767200601\\
47	1.52905767203099\\
48	1.52905767197249\\
49	1.52905767200774\\
50	1.52905767190004\\
};
\end{axis}

\end{tikzpicture}%
\captionsetup{justification=centering,font=scriptsize}
\caption{Comparison of the Frobenius norm approximation error for the noisy low-rank matrix. No power method, $q=0$, (left), and $q=2$ (right).} 
\label{fig:FroErr_1}       
\end{center}
\end{figure}

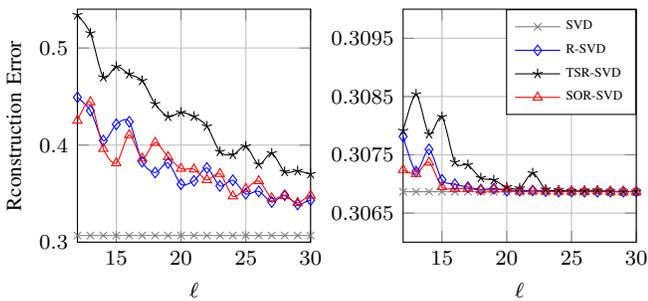
\begin{figure}[t]
\begin{center}
%
%
%
\usetikzlibrary{positioning,calc}

\definecolor{mycolor1}{rgb}{0.00000,1.00000,1.00000}%
\definecolor{mycolor2}{rgb}{1.00000,0.00000,1.00000}%

\pgfplotsset{every axis label/.append style={font=\footnotesize},
every tick label/.append style={font=\footnotesize}
}

\begin{tikzpicture}[font=\footnotesize] 

\begin{axis}[%
name=ber,
width  = 0.35\columnwidth,
height = 0.35\columnwidth,
scale only axis,
xmin  = 12,
xmax  = 30,
xlabel= {$\ell$},
xmajorgrids,
ymin = 0.3,
ymax = 0.54,
ylabel={Rconstruction Error},
ymajorgrids,
]


\addplot+[color=gray, mark=x]
table[row sep=crcr]{
12	0.306866152442753\\
13	0.306866152442753\\
14	0.306866152442753\\
15	0.306866152442753\\
16	0.306866152442753\\
17	0.306866152442753\\
18	0.306866152442753\\
19	0.306866152442753\\
20	0.306866152442753\\
21	0.306866152442753\\
22	0.306866152442753\\
23	0.306866152442753\\
24	0.306866152442753\\
25	0.306866152442753\\
26	0.306866152442753\\
27	0.306866152442753\\
28	0.306866152442753\\
29	0.306866152442753\\
30	0.306866152442753\\
31	0.306866152442753\\
32	0.306866152442753\\
33	0.306866152442753\\
34	0.306866152442753\\
35	0.306866152442753\\
36	0.306866152442753\\
37	0.306866152442753\\
38	0.306866152442753\\
39	0.306866152442753\\
40	0.306866152442753\\
}; 

\addplot+[smooth,color = blue,solid, every mark/.append style={solid}, mark=diamond]
  table[row sep=crcr]{
12	0.449205544606831\\
13	0.435437574589484\\
14	0.404631540478506\\
15	0.421360724115266\\
16	0.424090139936447\\
17	0.382843307464774\\
18	0.371744220751763\\
19	0.381513305651470\\
20	0.359675829380473\\
21	0.363194294552628\\
22	0.376435641300771\\
23	0.357966237675757\\
24	0.363769906314874\\
25	0.349862307538052\\
26	0.352354499065089\\
27	0.341334552629257\\
28	0.348430255468696\\
29	0.338816339004226\\
30	0.343249382176485\\
31	0.340224687578436\\
32	0.342261850149639\\
33	0.335632132367474\\
34	0.333001289895991\\
35	0.333163627156585\\
36	0.333382405800524\\
37	0.331896742011011\\
38	0.331294260450118\\
39	0.331806559899167\\
40	0.326173699434886\\
};

\addplot+[smooth,color=black,solid, every mark/.append style={solid}, mark=star]
  table[row sep=crcr]{
12	0.534048905202032\\
13	0.515224707938049\\
14	0.469934715267566\\
15	0.480817067088222\\
16	0.473106448732004\\
17	0.466272490344923\\
18	0.442409380895000\\
19	0.429348603919686\\
20	0.433631653077700\\
21	0.429459291185556\\
22	0.419297801019297\\
23	0.393198916450334\\
24	0.390137811108593\\
25	0.398714395640027\\
26	0.380103739060776\\
27	0.391654464739988\\
28	0.372586825587750\\
29	0.373386847218662\\
30	0.370261904781396\\
31	0.370378905770177\\
32	0.361767738136194\\
33	0.374887724754079\\
34	0.355479095964872\\
35	0.360102938234415\\
36	0.354243326072056\\
37	0.354125765238859\\
38	0.349318430875914\\
39	0.354282115721354\\
40	0.348995381549424\\
};

\addplot+[smooth,color=red,solid, every mark/.append style={solid}, mark=triangle]
  table[row sep=crcr]{
12	0.425344242346031\\
13	0.444361912880474\\
14	0.396305368505801\\
15	0.381490807108192\\
16	0.410484885305019\\
17	0.386260856687830\\
18	0.402531525554318\\
19	0.387829444663677\\
20	0.375827087844120\\
21	0.375131713450317\\
22	0.364187299691263\\
23	0.370470797822650\\
24	0.347389560583645\\
25	0.355049467704961\\
26	0.362986200164210\\
27	0.345364567882411\\
28	0.347804254265961\\
29	0.340505975590754\\
30	0.348108051278137\\
31	0.339205772996120\\
32	0.331270594382885\\
33	0.338784160180907\\
34	0.337946750677884\\
35	0.332888114895528\\
36	0.338027262185946\\
37	0.334240962722358\\
38	0.331714654313962\\
39	0.325277982731006\\
40	0.325843084427713\\
}; 
\end{axis}

\begin{axis}[%
name=SumRate,
at={($(ber.east)+(35,0em)$)},
		anchor= west,
width  = 0.35\columnwidth,
height = 0.35\columnwidth,
scale only axis,
xmin  = 12,
xmax  = 30,
xlabel= {$\ell$},
xmajorgrids,
ymin = 0.306,
ymax = 0.310,
ylabel={},
ymajorgrids,
ytick       ={0.3095, 0.3085, 0.3075, 0.3065 },
yticklabels ={$0.3095$, $0.3085$ , $0.3075$, $0.3065$},
legend entries={SVD, R-SVD, TSR-SVD,SOR-SVD},
legend style={at={(1,1)},anchor=north east,draw=black,fill=white,legend cell align=left,font=\tiny}
]

\addplot+[color=gray, mark=x]
table[row sep=crcr]{
12	0.306866152442753\\
13	0.306866152442753\\
14	0.306866152442753\\
15	0.306866152442753\\
16	0.306866152442753\\
17	0.306866152442753\\
18	0.306866152442753\\
19	0.306866152442753\\
20	0.306866152442753\\
21	0.306866152442753\\
22	0.306866152442753\\
23	0.306866152442753\\
24	0.306866152442753\\
25	0.306866152442753\\
26	0.306866152442753\\
27	0.306866152442753\\
28	0.306866152442753\\
29	0.306866152442753\\
30	0.306866152442753\\
31	0.306866152442753\\
32	0.306866152442753\\
33	0.306866152442753\\
34	0.306866152442753\\
35	0.306866152442753\\
36	0.306866152442753\\
37	0.306866152442753\\
38	0.306866152442753\\
39	0.306866152442753\\
40	0.306866152442753\\
};

\addplot+[smooth,color=blue,solid, every mark/.append style={solid}, mark=diamond]
  table[row sep=crcr]{
12	0.307812804457339\\
13	0.307212290971861\\
14	0.307593951876567\\
15	0.307069965134750\\
16	0.306993125697704\\
17	0.306945407631565\\
18	0.306904117334126\\
19	0.306921262928803\\
20	0.306884423899897\\
21	0.306882140176638\\
22	0.306887556183478\\
23	0.306877285746721\\
24	0.306868965677699\\
25	0.306870398259579\\
26	0.306866952639779\\
27	0.306867554197703\\
28	0.306866819722270\\
29	0.306866952139780\\
30	0.306867389948175\\
31	0.306866600504781\\
32	0.306866533676469\\
33	0.306866249631790\\
34	0.306866281678244\\
35	0.306866224997293\\
36	0.306866197935284\\
37	0.306866206285909\\
38	0.306866183440977\\
39	0.306866171897562\\
40	0.306866200591985\\
};

\addplot+[smooth,color=black, mark=star]
  table[row sep=crcr]{
12	0.307913907454331\\
13	0.308539451135885\\
14	0.307855312771677\\
15	0.308151067220466\\
16	0.307372680977667\\
17	0.307325846820717\\
18	0.307102444408104\\
19	0.307067016372668\\
20	0.306944528407247\\
21	0.306926151493813\\
22	0.307186242886879\\
23	0.306905022863745\\
24	0.306891485181461\\
25	0.306877648643810\\
26	0.306885430526996\\
27	0.306889258285182\\
28	0.306871430452282\\
29	0.306868603491994\\
30	0.306869983425075\\
31	0.306869006891054\\
32	0.306869545250910\\
33	0.306867283700062\\
34	0.306869259052593\\
35	0.306869666602328\\
36	0.306867863116743\\
37	0.306866901976025\\
38	0.306866466016204\\
39	0.306866588216440\\
40	0.306866916551265\\
};

\addplot+[smooth,color=red,solid, every mark/.append style={solid}, mark=triangle]
  table[row sep=crcr]{
12	0.307242945140612\\
13	0.307173868210809\\
14	0.307370372381765\\
15	0.306948874407070\\
16	0.306920667191400\\
17	0.306922084514712\\
18	0.306891216020234\\
19	0.306905723316579\\
20	0.306910880818569\\
21	0.306872790024743\\
22	0.306876228492318\\
23	0.306871937568673\\
24	0.306870751330206\\
25	0.306870056823430\\
26	0.306867888596211\\
27	0.306875875916978\\
28	0.306866826948342\\
29	0.306867149314994\\
30	0.306867179090114\\
31	0.306866599592673\\
32	0.306866413800793\\
33	0.306866244381133\\
34	0.306866238840110\\
35	0.306866207510223\\
36	0.306866176527391\\
37	0.306866324510101\\
38	0.306866244157585\\
39	0.306866185060986\\
40	0.306866173529286\\
};

\end{axis}

\end{tikzpicture}%
\captionsetup{justification=centering,font=scriptsize}
\caption{Comparison of the Frobenius norm approximation error for the matrix with polynomially decaying singular values. No power method, $q=0$, (left), and $q=2$ (right).} 
\label{fig:FroErr_2}       
\end{center}
\end{figure}

Figures \ref{fig:FroErr_1} and \ref{fig:FroErr_2} illustrate that when $q=0$, 
as $\ell$ increases, the approximation error $e_k$ declines for all algorithms. 
In this case, it is observed that the R-SVD and SOR-SVD algorithms have similar performance, exceeding the performance of the TSR-SVD algorithm. By incorporating two steps of the power iteration, however, at an increased computational cost, the accuracy of the low-rank approximations substantially improves for all algorithms. 

\subsection{Empirical Evaluation of SOR-SVD Error Bounds}
\label{subEmpEva}
\begin{figure}[t]
\begin{center}
%
%
%
\usetikzlibrary{positioning,calc}

\definecolor{mycolor1}{rgb}{0.00000,1.00000,1.00000}%
\definecolor{mycolor2}{rgb}{1.00000,0.00000,1.00000}%

\pgfplotsset{every axis label/.append style={font=\footnotesize},
every tick label/.append style={font=\footnotesize}
}

\begin{tikzpicture}[font=\footnotesize] 

\begin{axis}[%
name=ber,
width  = 0.35\columnwidth,
height = 0.35\columnwidth,
scale only axis,
xmin  = 21,
xmax  = 40,
xlabel= {$\ell$},
xmajorgrids,
ymin = 1.9,
ymax = 10.3,
ylabel={Relative Error},
ymajorgrids,
]


\addplot+[smooth,color = black,solid, every mark/.append style={solid}, mark=x]
  table[row sep=crcr]{
21	10.2973740767222\\
22	8.11246439413894\\
23	7.42707253235115\\
24	8.24848923047241\\
25	7.03561248742314\\
26	6.81548108768527\\
27	6.96311125838536\\
28	7.46107637902830\\
29	6.91480724810866\\
30	6.73462360437058\\
31	6.61371081159226\\
32	6.96484470526954\\
33	6.41239843509274\\
34	6.34148552620701\\
35	5.88988559072396\\
36	6.46697982813915\\
37	6.10190296297458\\
38	6.32558345866482\\
39	5.71041031520734\\
40	6.03565802732799\\
};

\addplot+[smooth,color=red,solid, every mark/.append style={solid}, mark=triangle]
  table[row sep=crcr]{
21	2.71969809316296\\
22	2.62724025605360\\
23	2.52823048648176\\
24	2.57458882646476\\
25	2.45089664684020\\
26	2.27509480633477\\
27	2.27931316514451\\
28	2.38206925806218\\
29	2.27576812920238\\
30	2.27722692116676\\
31	2.23184280478005\\
32	2.32195178342830\\
33	2.10073341267409\\
34	2.12769710890387\\
35	2.04203253089069\\
36	2.11788557375616\\
37	2.10008164366605\\
38	2.02124874224937\\
39	1.98557536185156\\
40	2.02591667783624\\
};

\end{axis}

\begin{axis}[%
name=SumRate,
at={($(ber.east)+(35,0em)$)},
		anchor= west,
width  = 0.35\columnwidth,
height = 0.35\columnwidth,
scale only axis,
xmin  = 21,
xmax  = 40,
xlabel= {$\ell$},
xmajorgrids,
ymin = 1.5314520,
ymax = 1.5314635,
ylabel={},
ymajorgrids,
ytick       ={1.531461, 1.531457 , 1.531453},
yticklabels ={$1.531461$, $1.531457$ , $1.531453$},
legend entries={Theoretical bound of Thr 5, 
$\|{\bf A} - \hat{\bf A}_{\text{SOR}}\|_{F}$},
legend style={at={(1,0.3)},anchor=north east,draw=black,fill=white,legend cell align=left,font=\tiny}
]


\addplot+[smooth,color = black,solid, every mark/.append style={solid}, mark=x]
  table[row sep=crcr]{
21	1.53146016657814\\
22	1.53146000549539\\
23	1.53145998221198\\
24	1.53146002322358\\
25	1.53146001663551\\
26	1.53145994976274\\
27	1.53145995523226\\
28	1.53145993519673\\
29	1.53145991939548\\
30	1.53145990756339\\
31	1.53145990988361\\
32	1.53145992376882\\
33	1.53145990164411\\
34	1.53145989147099\\
35	1.53145989397370\\
36	1.53145989732210\\
37	1.53145989577100\\
38	1.53145988784944\\
39	1.53145988552825\\
40	1.53145988988922\\
};

\addplot+[smooth,color=red,solid, every mark/.append style={solid}, mark=triangle]
  table[row sep=crcr]{
21	1.53145985619605\\
22	1.53145985103603\\
23	1.53145985058079\\
24	1.53145985142121\\
25	1.53145985087982\\
26	1.53145984965702\\
27	1.53145984963578\\
28	1.53145984955730\\
29	1.53145984924536\\
30	1.53145984896791\\
31	1.53145984875800\\
32	1.53145984895096\\
33	1.53145984844860\\
34	1.53145984815883\\
35	1.53145984825178\\
36	1.53145984834727\\
37	1.53145984830206\\
38	1.53145984808869\\
39	1.53145984810789\\
40	1.53145984805787\\
};

\end{axis}

\end{tikzpicture}%
\captionsetup{justification=centering,font=scriptsize}
\caption{Comparison of the Frobenius norm error of the SOR-SVD algorithm with
the theoretical bound (Theorem 5). No power method, $q=0$, (left), and $q=2$ (right).} 
\label{fig:ThrBo_F}       
\end{center}
\end{figure}
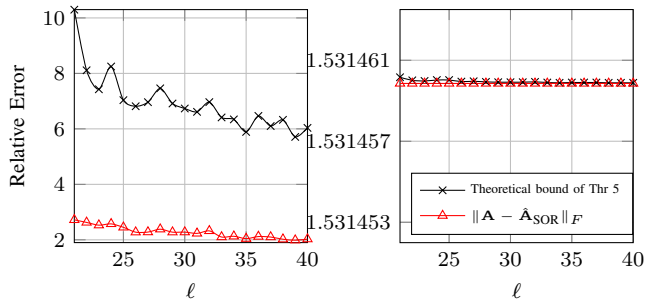

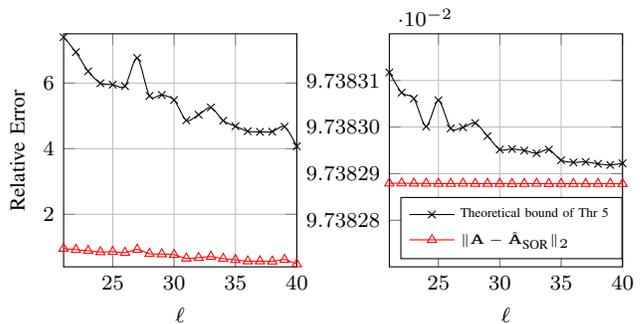
\begin{figure}[t]
\begin{center}
%
%
%
\usetikzlibrary{positioning,calc}

\definecolor{mycolor1}{rgb}{0.00000,1.00000,1.00000}%
\definecolor{mycolor2}{rgb}{1.00000,0.00000,1.00000}%

\pgfplotsset{every axis label/.append style={font=\footnotesize},
every tick label/.append style={font=\footnotesize}
}

\begin{tikzpicture}[font=\footnotesize] 

\begin{axis}[%
name=ber,
width  = 0.35\columnwidth,
height = 0.35\columnwidth,
scale only axis,
xmin  = 21,
xmax  = 40,
xlabel= {$\ell$},
xmajorgrids,
ymin = 0.4,
ymax = 7.5,
ylabel={Relative Error},
ymajorgrids,
]


\addplot+[smooth,color = black,solid, every mark/.append style={solid}, mark=x]
  table[row sep=crcr]{
21	7.39920667088094\\
22	6.94222714561974\\
23	6.35793695157030\\
24	5.99521915894009\\
25	5.95245401941743\\
26	5.90851369549741\\
27	6.76381949884114\\
28	5.60837588616217\\
29	5.63791413453956\\
30	5.48375124142166\\
31	4.86411319973565\\
32	5.03475682763435\\
33	5.25492321499383\\
34	4.85397346943800\\
35	4.68374983199101\\
36	4.53211619990813\\
37	4.51281100431377\\
38	4.51841741933881\\
39	4.67094452219417\\
40	4.07623420660115\\
};

\addplot+[smooth,color=red,solid, every mark/.append style={solid}, mark=triangle]
  table[row sep=crcr]{
21	0.958802725019479\\
22	0.921725853540813\\
23	0.890225870556153\\
24	0.851831970134506\\
25	0.861225797228616\\
26	0.831662793149209\\
27	0.927413934701938\\
28	0.800429287606455\\
29	0.788010615665046\\
30	0.768483463366721\\
31	0.657649116312127\\
32	0.676733767046283\\
33	0.711156211734294\\
34	0.646400451434940\\
35	0.613843414944460\\
36	0.572312886014707\\
37	0.571926160390907\\
38	0.566458990492952\\
39	0.618306356663180\\
40	0.483483753208676\\
};

\end{axis}

\begin{axis}[%
name=SumRate,
at={($(ber.east)+(35,0em)$)},
		anchor= west,
width  = 0.35\columnwidth,
height = 0.35\columnwidth,
scale only axis,
xmin  = 21,
xmax  = 40,
xlabel= {$\ell$},
xmajorgrids,
ymin = 0.0973827,
ymax = 0.0973832,
ylabel={},
ymajorgrids,
ytick       ={0.0973831, 0.0973830 , 0.0973829, 0.0973828},
yticklabels ={$9.73831$, $9.73830$ , $9.73829$, $9.73828$},
legend entries={Theoretical bound of Thr 5, 
$\|{\bf A} - \hat{\bf A}_{\text{SOR}}\|_2$},
legend style={at={(1,0.3)},anchor=north east,draw=black,fill=white,legend cell align=left,font=\tiny}
]

\addplot+[smooth,color = black,solid, every mark/.append style={solid}, mark=x]
  table[row sep=crcr]{
21	0.0973831172731646\\
22	0.0973830739007262\\
23	0.0973830610025062\\
24	0.0973830011834551\\
25	0.0973830576135872\\
26	0.0973829971097660\\
27	0.0973829993499755\\
28	0.0973830086353606\\
29	0.0973829805905660\\
30	0.0973829517723971\\
31	0.0973829529361360\\
32	0.0973829494552501\\
33	0.0973829438401326\\
34	0.0973829517311263\\
35	0.0973829289704595\\
36	0.0973829241839130\\
37	0.0973829252892336\\
38	0.0973829210341949\\
39	0.0973829186725039\\
40	0.0973829222584503\\
};

\addplot+[smooth,color=red,solid, every mark/.append style={solid}, mark=triangle]
  table[row sep=crcr]{
21	0.0973828792582738\\
22	0.0973828794584715\\
23	0.0973828791134601\\
24	0.0973828790379216\\
25	0.0973828787081896\\
26	0.0973828785772500\\
27	0.0973828790589174\\
28	0.0973828786012769\\
29	0.0973828785885697\\
30	0.0973828785285217\\
31	0.0973828783158228\\
32	0.0973828784120796\\
33	0.0973828788005354\\
34	0.0973828782882339\\
35	0.0973828783875565\\
36	0.0973828784120977\\
37	0.0973828782865577\\
38	0.0973828784451357\\
39	0.0973828783856673\\
40	0.0973828782726668\\
};

\end{axis}

\end{tikzpicture}%
\captionsetup{justification=centering,font=scriptsize}
\caption{Comparison of the spectral norm error of the SOR-SVD algorithm with
the theoretical bound (Theorem 5). No power method, $q=0$, (left), and $q=2$ (right).} 
\label{fig:ThrBo_L2}       
\end{center}
\end{figure}

The theoretical error bounds for the SOR-SVD algorithm are given in Theorem \ref{Thr2}. To evaluate the accuracy of Theorem \ref{Thr2}, we form an input 
matrix according to \ref{NoisyLRM}. With the rank $k=20$ fixed, we increase 
the sample size parameter $\ell$, considering the assumption $2\le p\le \ell-k$. 
A comparison between the theoretical bounds and what are achieved in practice 
is shown in Figures \ref{fig:ThrBo_F} and \ref{fig:ThrBo_L2}; Figure \ref{fig:ThrBo_F} compares the Frobenius norm error with the corresponding theoretical bound, and Figure \ref{fig:ThrBo_L2} compares the spectral norm error with the corresponding theoretical bound.  

The effect of the power scheme can be easily seen from the figures; when $q=2$, 
the theoretical bounds given in Theorem \ref{Thr2} closely track the error in 
the low-rank approximation of Alg. \ref{Alg3}. We conclude that for the 
noisy low-rank matrix, the theoretical error bounds are empirically sharp. 
 
\subsection{Robust Principal Component Analysis}
\label{Exrpca}
Principal component analysis (PCA) \cite{Jackson91} is a linear dimensionality reduction technique that tranforms the data to a lower-dimensional 
subspace that captures the most features of the data. However, PCA is known 
to be very sensitive to grossly corrupted observations. After a long 
line of research to robustifying PCA against grossly perturbed observations, 
robust PCA \cite{WPMGR2009,CSPW2009, CLMW2009} was proposed. Robust PCA 
represents a corrupted low-rank matrix ${\bf X} \in \mathbb R^{m \times n}$ 
as a linear combination of a low-rank matrix ${\bf L}$ and a sparse matrix 
of corrupted entries ${\bf S}$ such that ${\bf X=L+S}$, by solving the 
following convex program:  
\begin{equation}
\begin{aligned}
&{\text{minimize}_{\bf(L, S)}} \ {\|{\bf L}\|_* + \lambda\|{\bf S}\|_1} \\
&{\text{subject to}} \ {\bf L+S=X},
\end{aligned}\label{equV1}
\end{equation} 
where, for any matrix $\bf B$, ${\|\mbox{\bf B}\|_*}  \triangleq \sum_i\sigma_i (\mbox{\bf B}) $ is the 
nuclear norm of the matrix $\mbox{\bf B}$ (sum of the singular values) and 
${\|\mbox{\bf B}\|_1} \triangleq \sum_{ij} |\mbox{\bf B}_{ij}|$ is the 
$\ell_{1}$-norm of the matrix  $\mbox{\bf B}$, and $\lambda>0$ is a weighting parameter. To solve \eqref{equV1} the method of augmented Lagrange multipliers 
(ALM) \cite{LLS2011, YY2009} iteratively minimizes the following augmented Lagrangian function with respect to either $\bf L$ or $\bf S$ with the other variable fixed to give the pair of optimal solutions $(\bf L^*, S^*)$:
\begin{equation}
\begin{aligned}
\mathcal{L}({\bf L}, {\bf S}, {\bf Y}, \mu) \triangleq 
& {\|{\bf L}\|_* + \lambda\|{\bf S}\|_1} + \langle\, {\bf Y}, {\bf X} - {\bf L}\ -{\bf S} \rangle \\
+ & \frac{\mu}{ 2}\|{\bf X - \  L\ -  S}\|_F^2,
\end{aligned} \label{equ4} 
\end{equation} where $ {\bf Y} \in \mathbb R^{m \times n}$ is a matrix of Lagrange multipliers, and $\mu>0$ is a penalty parameter.
To speed up the convergence of the ALM method, a continuation technique proposed in \cite{Toh2010} is used, i.e., in each iteration $\mu$ is increased such that $\mu_{k+1}=\rho\mu_k$, for some numerical constant $\rho>0$.

The bottleneck of the ALM method, however, is to compute an SVD at each iteration to approximate the low-rank component $\bf L$ of the input matrix $\bf X$, which is computationally demanding. Furthermore, as explained, the SVD techniques are challenging to be parallelized, and cannot easily be modified for partial factorizations \cite{Demmel97, GolubVanLoan96}. 

We therefore, by retaining the original objective function, apply the SOR-SVD as a surrogate to the SVD to solve the RPCA problem. The pseudocode of the proposed method, which we call 
ALM-SOR-SVD method, is given in Table \ref{TableThree}.  
\begin{table}[!htb]
\caption{Pseudo-code for RPCA solved by the ALM-SOR-SVD method.}
\normalsize
\label{TableThree} 
\noindent\hfil\rule{0.9\linewidth}{.2pt}\hfil

\begin{center}
\begin{minipage}{0.9\linewidth}
\renewcommand{\algorithmicrequire}{\textbf{Input:}} 
\begin{algorithmic}[1]
\REQUIRE ~~ 
Matrix ${\bf X}, \lambda, \mu_0,{\bar \mu}, \rho, {\bf Y}_0, {\bf S}_0, k=0$;
 \renewcommand{\algorithmicrequire}{\textbf{Output:}} 
 \REQUIRE ~~ 
 Low-rank plus sparse matrix
\WHILE {the algorithm does not converge}
        \STATE $({\bf U},{\bf Z},{\bf V}) = {\text{sor-svd}}({\bf X} 
        - {\bf S}_k +\mu_k^{-1} {\bf Y}_k)$;
        \STATE ${\bf L}_{k+1} = {\bf U} \mathcal{S}_{\mu_k^{-1}}
        ({\bf Z}){\bf V}^T$;
        \STATE ${\bf S}_{k+1} = \mathcal{S}_{\lambda\mu_k^{-1}}
        ({\bf X} - {\bf L}_{k+1} +\mu_k^{-1} {\bf Y})$;
        \STATE ${\bf Y}_{k+1} = {\bf Y}_k +\mu_k ({\bf X} - {\bf L}_{k+1} 
        - {\bf S}_{k+1})$;
        \STATE $\mu_{k+1} \leftarrow \text{max}(\rho\mu_k, {\bar \mu})$;
\ENDWHILE
\RETURN $\bf L^*$ and $\bf S^*$ 
\vspace*{-1.8em}
\end{algorithmic}
\end{minipage}
\end{center}

\begin{center}
\end{center}
\noindent\hfil\rule{0.9\linewidth}{.5pt}\hfil
\end{table}

Here $\mathcal{S}_\varepsilon (x) = {\text{sgn}(x)\text{max}}(|x| - \varepsilon, 0)$ is a soft-thresholding operator \cite{Hale2008}, and $ \lambda$, $\mu_0$, ${\bar \mu}$, $\rho$, ${\bf Y}_0$, and ${\bf S}_0$ are initial values.

In the next subsections, we verify the efficiency and efficacy of the ALM-SOR-SVD to solve the RPCA problem on randomly generated data, as well as real-time data. We compare the experimental results obtained with those of applying the partial SVD (by using PROPACK package) \cite{Larsen98}. The PROPACK function provides an efficient algorithm, suitable for approximating large low-rank matrices, which computes a specified number of largest singular values and corresponding singular 
vectors of a matrix by making use of the Lanczos bidiagonalization algorithm with partial reorthogonalization (BPRO). We run the experiments in MATLAB on a desktop PC with a 3 GHz intel Core i5-4430 processor and 8 GB of memory.\looseness-1

\subsubsection{Synthetic Data Recovery}
\label{subsecSyDataRec}

We generate a rank-$k$ matrix $\bf X = L + S$ as a sum of a low-rank matrix ${\bf L} \in \mathbb R^{n \times n}$ and a sparse matrix  ${\bf S}\in \mathbb R^{n \times n}$. The matrix ${\bf L}$ is generated by a matrix multiplication ${\bf L}={\bf U}{\bf V}^T$, where ${\bf U}$, ${\bf V} \in \mathbb R^{n \times k}$ have standard Gaussian distributed entries. The matrix ${\bf S}$ has $s$ non-zero entries independently drawn from the set $\lbrace$-50, 50$\rbrace$. We apply the ALM-SOR-SVD method, and the ALM method using the partial SVD, hereafter ALM-PSVD, on $\bf X$ to recover $\hat {\bf L}$ and $\hat {\bf S}$.\looseness-1 

Table \ref{TableFour} summarizes the numerical results where $\text{rank}
({\bf L})=0.05\times n$, and $s = \|{\bf S}\|_0=0.05\times n^2$, and Table \ref{TableFive} presents the results for a more challenging 
scenario where $\text{rank}({\bf L})= 0.05\times n$ and 
$s = \|{\bf S}\|_0=0.1\times n^2$. The results of ALM-PSVD method are reported 
in the numerators, and and those of the ALM-SOR-SVD method in the denominators.
In the Tables, $Time$ denotes the computational time in seconds, $Iter.$ 
denotes the number of iterations, and $\xi$ denotes the relative error. 
For the simulations, the initial values suggested in \cite{LLS2011} are adopted, and 
the algorithms stop when the relative error $\xi$ defined as ${\|{\bf X}-\hat{\bf L} - \hat{\bf S}\|_F}< 10^{-7} {\|{\bf X}\|_F}$ is satisfied. 

Since both SOR-SVD and the truncated SVD algorithms require a predetermined 
rank $\ell$ to perform the decomposition, we set $\ell=2r$, a random start, 
and $q=1$ for the SOR-SVD. 
\begin{table}[!htb]
\centering
\caption{\textnormal{Comparison of the ALM-SOR-SVD and ALM-PSVD methods 
for synthetic data recovery, for the case $r({\bf L})=0.05\times n$ and
 $s = 0.05\times n^2$.}}
\begin{tabular}{p{0.2cm} c p{0.3cm} c p{0.7cm} c c c}
\noindent\rule{7.9cm}{0.4pt}\\
$n$ & $r(\bf L)$ & $r(\hat{\bf L})$&$\|{\bf S}\|_0$&$\|\hat{\bf S}\|_0$&Time &Iter.&$\xi$ \\
\noindent\rule{7.9cm}{0.4pt}\\
500& 25 & $\dfrac{25}{25}$ & 12500 & $\dfrac{12500}{12500}$ & $\dfrac{2.2}{0.1}$ &
$\dfrac{17}{17}$ & $\dfrac{6.6\text{e}-8}{6.3\text{e}-8}$  \\
\\
1000& 50 & $\dfrac{50}{50}$ & 50000 & $\dfrac{50000}{50000}$ & $\dfrac{22}{0.7}$ &
$\dfrac{17}{17}$ & $\dfrac{5.4\text{e}-8}{5.3\text{e}-8}$  \\
\\
2000 & 100 & $\dfrac{100}{100}$ & 200000& $\dfrac{200000}{200000}$ & 
$\dfrac{167.9}{4.6}$ & 
$\dfrac{17}{17}$ &$\dfrac{5\text{e}-8}{5.1\text{e}-8}$ \\
\\
3000 & 150 & $\dfrac{150}{150}$ & 450000 & $\dfrac{450000}{450000}$ & 
$\dfrac{563.4}{11.8}$ &
$\dfrac{17}{17}$ & $\dfrac{4.8\text{e}-8}{4.9\text{e}-8}$ \\
\noindent\rule{7.9cm}{0.4pt} 
\end{tabular}
\label{TableFour}
\end{table}

\begin{table}[!htb]
\centering
\caption{\textnormal{Comparison of the ALM-SOR-SVD and ALM-PSVD methods for synthetic 
data recovery, for the case $r({\bf L})=0.05\times n$ and $s = 0.1\times n^2$.}}
\begin{tabular}{p{0.2cm} c p{0.3cm} c p{0.7cm} c c c}
\noindent\rule{7.9cm}{0.4pt}\\
$n$ & $r(\bf L)$ & $r(\hat{\bf L})$&$\|{\bf S}\|_0$&$\|\hat{\bf S}\|_0$&Time &Iter.&$\xi$ \\
\noindent\rule{7.9cm}{0.4pt}\\
500& 25 & $\dfrac{25}{25}$ & 25000 & $\dfrac{25000}{25000}$ & 
$\dfrac{2.5}{0.2}$ &
$\dfrac{20}{20}$ & $\dfrac{3.2\text{e}-8}{4.9\text{e}-8}$  \\
\\
1000& 50 & $\dfrac{50}{50}$ & 100000 & $\dfrac{100000}{100000}$ & 
$\dfrac{25.9}{0.8}$ & $\dfrac{19}{19}$ & $\dfrac{8.1\text{e}-8}{8.3\text{e}-8}$ \\
\\
2000 & 100 & $\dfrac{100}{100}$ & 400000& $\dfrac{400000}{400000}$ &  
$\dfrac{189}{5.3}$ 
& $\dfrac{19}{19}$& $\dfrac{6.8\text{e}-8}{7.4\text{e}-8}$  \\
\\
3000 & 150 & $\dfrac{150}{150}$ & 900000 & $\dfrac{900000}{900000}$ & 
$\dfrac{609.2}{13.6}$ 
& $\dfrac{19}{19}$ & $\dfrac{7.2\text{e}-8}{7.6\text{e}-8}$ \\
\noindent\rule{7.9cm}{0.4pt}\\
\end{tabular}
\label{TableFive}
\end{table}

The results in Tables \ref{TableFour} and \ref{TableFive} lead us 
to make several conclusions on the ALM-SOR-SVD method:
\begin{itemize}
\item it successfully detects the numerical rank $k$ in all cases.
\item it provides the exact recovery of the sparse matrix $\bf S$ from 
$\bf X$, with the same number of iterations compared to the ALM-PSVD method. 
\item In terms of runtime, it outperforms the ALM-PSVD method with 
speedups of up to $47$ times.
\end{itemize}

\subsubsection{Background Subtraction in Surveillance Video}
\label{subsecBSub}
Separating the foreground or moving objects from the background in a video 
sequence fits nicely into the RPCA model, where the background is modeled 
by a low-rank matrix, and moving objects are modeled by a sparse matrix.

In our experiment, we use four different real-time videos introduced in \cite{LHGT2004}. 
The first video sequence has 200 grayscale frames with dimensions ${176 \times 144}$ 
in each frame, and has been taken in a hall of an airport. We stack each frame as 
a column of the data matrix ${\bf X} \in \mathbb R^{25344 \times 200}$. The second 
video has 200 grayscale frames with dimensions ${256 \times 320}$ in each frame, 
and has been taken in a shopping mall. Thus ${\bf X} \in \mathbb R^{81920 \times 200}$. These two videos have relatively static background. The third video has 200 
grayscale frames with dimensions ${130 \times 160}$ in each frame, i.e., ${\bf X} \in \mathbb R^{20800 \times 200}$, and has been taken from an escalator at an airport. The background of this video changes due to the moving escalator. The fourth video has 250 grayscale frames with dimensions ${128 \times 160}$ in each frame taken in an office. Thus ${\bf X} \in \mathbb R^{20480 \times 250}$. 
In this video, the illumination changes drastically, making it very challenging
to analyze.

The predetermined rank $\ell$, assigned to both algorithms, is obtained by invoking the following bound that, for any rank-$k$ matrix $\bf B$, relates $k$ with the Frobenius and nuclear norms \cite{GolubVanLoan96}:\looseness-1  
\begin{equation}
\|{\bf B}\|_* \le \sqrt{k}\|{\bf B}\|_F.
\label{equ21}
\end{equation}

We assign the value of $k$ to $\ell$, i.e., $\ell=k$, and  set $q=1$ for 
the SOR-SVD algorithm.

\begin{figure}[!htb]
\centering
\includegraphics[width=0.48\textwidth,height=4cm]{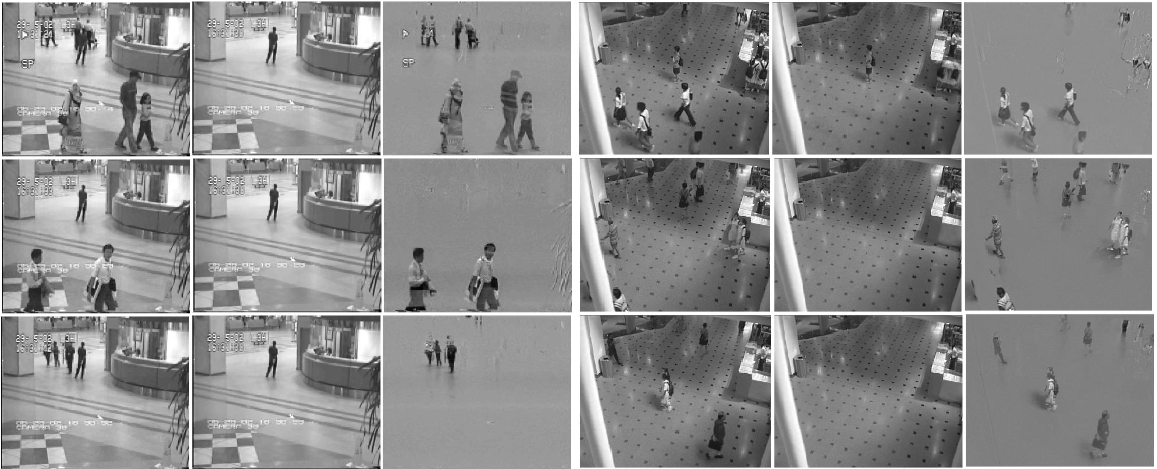}
\caption{Images in columns 1 and 4 are frames of the video sequence of an 
airport and a shopping mall, respectively. Images in columns 2 and 5  are 
recovered backgrounds $\hat {\bf L}$, and columns 3 and 6 correspond 
to foregrounds $\hat {\bf S}$.}
\label{fig_AirSh}
\end{figure}

\begin{figure}[!htb]
\centering
\includegraphics[width=0.48\textwidth,height=4cm]{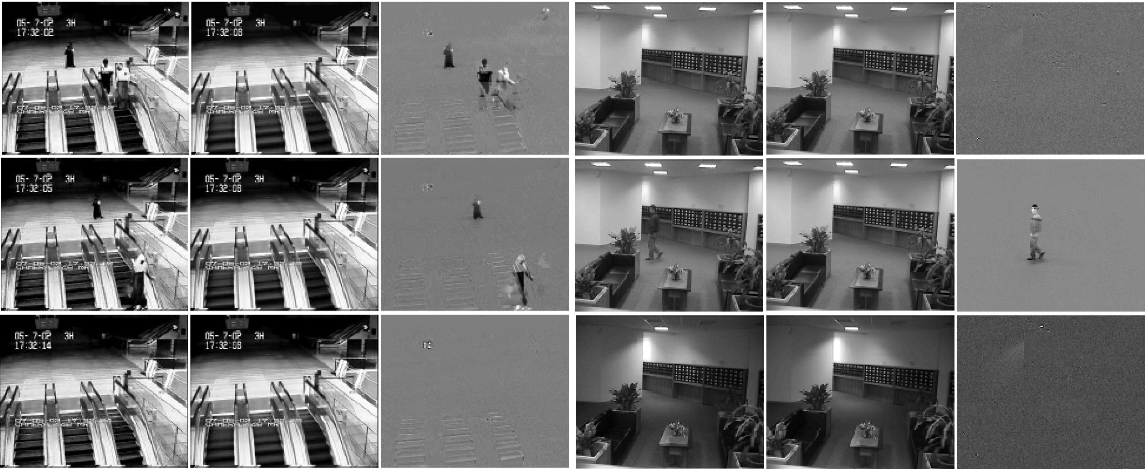}
\caption{Images in columns 1 and 4 are frames of the video sequence of an 
escalator and an office, respectively. Images in columns 2 and 5  are 
recovered backgrounds $\hat {\bf L}$, and columns 3 and 6 correspond 
to foregrounds $\hat {\bf S}$.}
\label{fig_EscOff}
\end{figure}

Some sample frames of the videos with corresponding recovered backgrounds and foregrounds are shown in Figures \ref{fig_AirSh} and \ref{fig_EscOff}. 
We only show the results of the ALM-SOR-SVD method since the results returned 
by the ALM-PSVD method are visually identical. It is evident that the proposed method can successfully recover the low-rank and sparse components of the data matrix in all scenarios.

Table \ref{TableSix} summarizes the numerical results. In all cases, 
the ALM-SOR-SVD method outperforms the ALM-PSVD method in terms of runtime, 
while having the same number of iterations.

\begin{table}[!htb]
\centering
\caption{\textnormal{Comparison of the ALM-SOR-SVD and ALM-PSVD methods 
for real-time data recovery.}}
\begin{tabular}{p{1.8cm} p{2cm} p{0.4cm}cc} 
\noindent\rule{7.5cm}{0.4pt}\\
Dataset & Dimensions &Time & Iter. &  $\xi$  \\
\noindent\rule{7.5cm}{0.4pt}\\
Airport hall & $25344\times 200 $ & $\dfrac{14.2}{5.7}$ & $\dfrac{36}{36}$
& $\dfrac{6\text{e}-8}{6.6\text{e}-8}$  \\
\\
Shopping mall & $81920 \times 200 $ & $\dfrac{44.2}{19.1}$ & $\dfrac{19}{19}$ &
$\dfrac{6.9\text{e}-8}{6.8\text{e}-8}$ \\
\\
Escalator & $20800\times 200 $ & $\dfrac{11.3}{4.6}$ & $\dfrac{36}{36}$ &
$\dfrac{7.5\text{e}-8}{6.5\text{e}-8}$ \\
\\
Lobby & $20480 \times 250$  & $\dfrac{11.2}{5.7}$  & $\dfrac{36}{36}$ 
& $\dfrac{6\text{e}-8}{6.6\text{e}-8}$\\
\\
Yale B03 & $32256 \times 64$ & $\dfrac{6}{2.5}$ & $\dfrac{36}{36}$
& $\dfrac{9.2\text{e}-8}{7.5\text{e}-8}$ \\
\\
Yale B08 & $32256 \times 64$ & $\dfrac{7.2}{2.6}$ & $\dfrac{36}{36}$
& $\dfrac{9.8\text{e}-8}{7.6\text{e}-8}$ \\
\noindent\rule{7.5cm}{0.4pt}\\
\end{tabular}
\label{TableSix}
\end{table}

\subsubsection{Shadow and Specularity Removal From Face Images}
\label{subsecShRemoval}
Another task in computer vision that can be modeled as a robust PCA problem is removing shadows and specularities from face images; there exist images 
of the same face taken under varying illumination, where the face image can be modeled by a low-rank component and shadows and specularities can be modeled by a sparse component. 

\begin{figure}[!htb]
\centering
\includegraphics[width=0.45\textwidth,height=4cm]{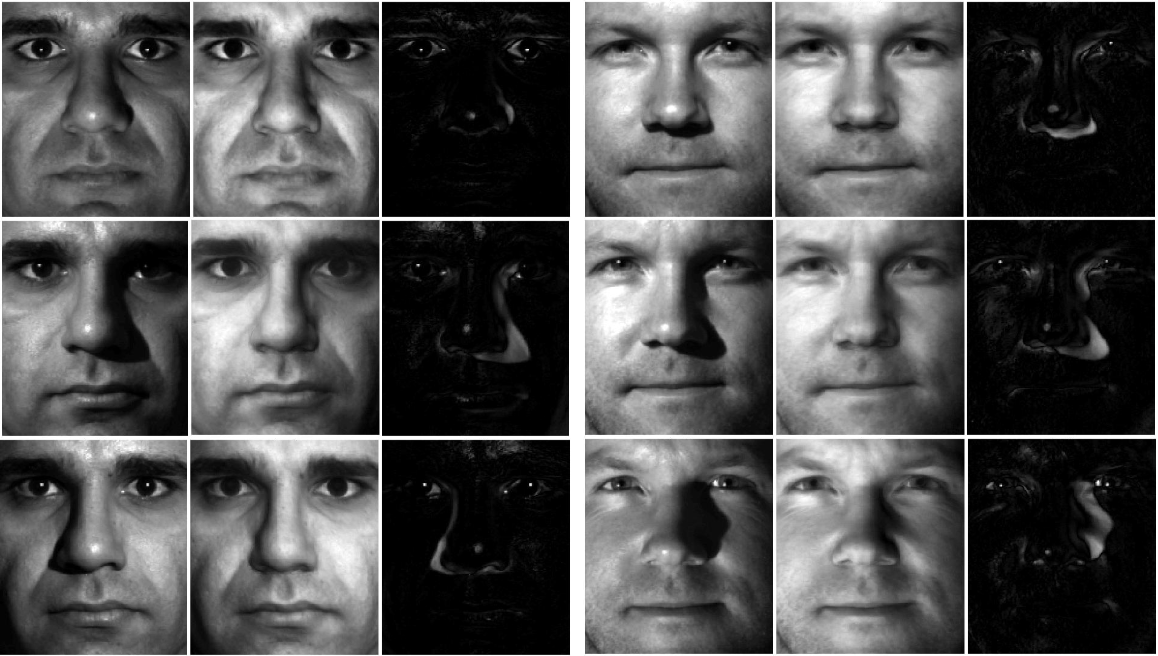}
\caption{Images in columns 1 and 4 are face images under different illuminations. Images in columns 2 and 5 are are recovered images after removing shadows and specularities, and images in columns 3 and 6 correspond to the removed shadows 
and specularities.}
\label{figFace}
\end{figure}

In our experiment, we use face images taken from the Yale B face database \cite{Georghiades2001}. Each image has dimensions ${192 \times 168}$ 
with a total of 64 illuminations. The images are stacked as columns of the data matrix ${\bf X} \in \mathbb R^{32256 \times 64}$. The recovered images are shown in Figure \ref{figFace}, and the numerical results are presented in Table \ref{TableSix}.

We conclude that the ALM-SOR-SVD method can successfully recovers the 
low-rank and sparse matrices from the dataset with speedups of up to 
$2.7$ times, compared to the ALM-PSVD method with the same number of iterations.

\section{Conclusion}
\label{secCon}

In this paper we have proposed a randomized algorithm termed SOR-SVD to 
compute a low-rank approximation of a given matrix. We have presented an 
error analysis of the SOR-SVD algorithm, and experimentally verified that 
the error bound is sharp for one type of low-rank matrices, whereas our 
empirical evaluations demonstrate the effectiveness of the SOR-SVD for several 
types of low-rank matrices.

For a given matrix, the performance of the SOR-SVD algorithm exceeds the 
performance of the TSR-SVD algorithm, while being similar to that of the R-SVD. However, the SOR-SVD can exploit modern computational platforms better by exposing higher levels of parallelism than the R-SVD.

We also applied the SOR-SVD algorithm to solve the RPCA problem via the ALM 
method. Our empirical studies show that the resulting ALM-SOR-SVD method 
and ALM-PSVD, obtained by applying the partial SVD in the ALM method, have 
similar behavior, while the ALM-SOR-SVD is substantially faster. 

\appendices
\section{Proof of Proposition \ref{Prop_G}}
\label{secAProp_G}
To prove Proposition \ref{Prop_G}, we first present several key results that are used later on. 
 
\begin{proposition}(Halko et al. \cite{HMT2009}). For fixed matrices ${\bf S}$ and
${\bf T}$, and a standard Gaussian matrix ${\bf G}$, we have
\begin{equation}
\mathbb{E}\|{\bf SGT}\|_2 \le \|{\bf S}\|_2\|{\bf T}\|_F+\|{\bf S}\|_F\|{\bf T}\|_2.
\notag
\end{equation}
\label{PropHalkoG} 
\end{proposition}

\begin{proposition}(Halko et al. \cite{HMT2009}). Let $h$ be a real valued Lipschitz function on matrices: 
\begin{equation}
|h({\bf X}) - h({\bf X})| \le L\|{\bf X}-{\bf Y}\|_F, \quad \text{for all} \quad 
{\bf X}, {\bf Y} 
\notag
\end{equation}
where $L>0$. Draw a standard Gaussian matrix $\bf G$. Then
\begin{equation}
\mathbb{P} \{ h({\bf G}) \ge \mathbb{E} h({\bf G}) +Lu\} \le \text{e}^{-u^2/2}.
\notag
\end{equation}
\label{PropHalkoL}
\end{proposition}

\begin{proposition}(Halko et al. \cite{HMT2009}). Let ${\bf G}\in 
\mathbb R^{\ell-p \times p}$ be a Gaussian matrix, where $p \ge 0$ 
and $\ell-p \ge 2$. Then for $t \ge 1$, 
\begin{equation}
 \mathbb{P}\Big\{\|{\bf G}^{\dagger}\|_2 \ge 
 \dfrac{\text{e}t\sqrt{\ell}}{p+1}\Big\} \le t^{-(p+1)}.
 \notag
\end{equation}
\label{PropHalkoP-I}
\end{proposition}

\begin{proposition}(Gu \cite{Gu2015}). Let $g(\cdot)$ be a non-negative continuously 
differentiable function with $g(0) = 0$, and let $\bf G$ be a random matrix, then 
\begin{equation}
\mathbb{E}_g(\| {\bf G}\|_2) = \int_{0}^{\infty} 
g\prime(x)\mathbb{P}\{\|{\bf G}\|_2\ \ge x\}dx.
\notag
\end{equation}
\label{Gu_g}
\end{proposition}

We first define the following function  
\begin{equation}
\begin{aligned}
g(x) & \triangleq \sqrt{\dfrac{\alpha^2 x^2}{1 + \beta^2x^2}},
\end{aligned}
\notag
\end{equation}
whose derivative is defined as
\begin{equation}
\begin{aligned}
g\prime(x) & = \dfrac{\alpha^2x}{(1 + \beta^2x^2)^2\sqrt{\dfrac{\alpha^2 x^2}{1 + \beta^2x^2}}}.
\end{aligned}
\notag
\end{equation}
where $\alpha, \beta>0$.

Next, for the Gaussian matrix ${\bf G}\in \mathbb R^{m \times n}$, we 
define a function $h({\bf G})=\|{\bf G}\|_2$. By Proposition \ref{PropHalkoG}, 
it follows that
\begin{equation}
\mathbb{E}(h({\bf G})) \le \sqrt{m}+\sqrt{n}< \sqrt{m}+\sqrt{n}+3 
\triangleq \varepsilon.
\notag
\end{equation}

By definition of $g(x)$, Proposition \ref{Gu_g}, and Proposition 
\ref{PropHalkoL}, for $x = u - \varepsilon$, we have 
\begin{equation}
\mathbb{P}\{\|{\bf G}\|_2 \ge x\} \le \text{e}^{-u^2/2}.
\notag
\end{equation}

We can rewrite equation \eqref{eqProp_G} as follows 
\begin{equation}
\begin{aligned}
& \mathbb{E}\Bigg(\sqrt{\dfrac{\alpha^2 \|{\bf G}\|_2^2}{1 + \beta^2\|{\bf G}\|_2^2}}\Bigg) = \mathbb{E}(g(\|{\bf G}\|_2)) \\ & =\int_{0}^{\infty} g\prime(x)\mathbb{P}\{\|{\bf G}\|_2 \ge x\}dx \\
& \le \int_{0}^{\varepsilon} g\prime(x)dx +
\int_{\varepsilon}^{\infty} g\prime(x)\mathbb{P}\{\|{\bf G}\|_2 \ge x\}dx \\
& \le \sqrt{\dfrac{\alpha^2 \varepsilon^2}{1 + \beta^2\varepsilon^2}} + 
\int_{\varepsilon}^{\infty} \dfrac{\alpha^2x}{(1 + \beta^2x^2)^2\sqrt{\dfrac{\alpha^2 x^2}{1 + \beta^2x^2}}}\text{e}^{({x-\varepsilon})^2/2}dx\\
& = \sqrt{\dfrac{\alpha^2 \varepsilon^2}{1 + \beta^2\varepsilon^2}} + 
\dfrac{\alpha^2}{(1 + \beta^2\varepsilon^2)^2\sqrt{\dfrac{\alpha^2 \varepsilon^2}{1 + \beta^2\varepsilon^2}}} \times \\ 
& \quad \underbrace{\int_{0}^{\infty}(u + \varepsilon)^{{-u}^2/2}du}_{\varepsilon\sqrt{\pi/2}+1}.
\end{aligned}
\notag
\end{equation}

We now must find a $\nu>0$ such that 

\begin{equation}
\begin{aligned}
\sqrt{\dfrac{\alpha^2 \varepsilon^2}{1 + \beta^2\varepsilon^2}} + 
\dfrac{\alpha^2(\varepsilon\sqrt{\pi/2}+1)}{(1 + \beta^2\varepsilon^2)^2\sqrt{\dfrac{\alpha^2 \varepsilon^2}{1 + \beta^2\varepsilon^2}}} \le 
\sqrt{\dfrac{\alpha^2 \nu^2}{1 + \beta^2 \nu^2}},
\end{aligned}
\notag
\end{equation}
which leads to
\begin{equation}
\begin{aligned} 
\nu^2-\varepsilon^2 \ge &  \dfrac{1 + \beta^2 \nu^2}{1 + \beta^2 \varepsilon^2}
\times (\varepsilon\sqrt{\pi/2}+1) \times 
\left[ \dfrac{\sqrt{\dfrac{\alpha^2 \nu^2}{1 + \beta^2 \nu^2}}}
{\sqrt{\dfrac{\alpha^2 \varepsilon^2}{1 + \beta^2 \varepsilon^2}}} +1 \right].
\end{aligned}
\notag
\end{equation}

The right-hand side of the inequality approaches the maximum value as $\beta$ approaches $\infty$. Thus $\nu$ must satisfy

\begin{equation}
\begin{aligned}
\nu^2-\varepsilon^2 \ge \dfrac{\nu^2}{\varepsilon^2}(\varepsilon\sqrt{\pi/2}+1).
\end{aligned}
\notag
\end{equation}

which results in
\begin{equation}
\begin{aligned}
\nu\ge \dfrac{\varepsilon^2}{\sqrt{\varepsilon^2 - (\varepsilon\sqrt{\pi/2}+1)}}.
\end{aligned}
\notag
\end{equation}

The inequality is satisfied when $\nu = \sqrt{m}+\sqrt{n}+7 = \varepsilon +4$. 

\section{Proof of Proposition \ref{Prop_PI}}
\label{secAProp_PI}
According to Proposition \ref{PropHalkoP-I}, for any $x>0$, we have

\begin{equation}
 \mathbb{P}\Big\{\|{\bf G}^{\dagger}\|_2 \ge x\} \le 
 \Big(\dfrac{p+1}{\text{e}\sqrt{\ell}}x\Big)^{-(p+1)}.
 \notag
\end{equation}

With similar arguments to those in the proof of Proposition \ref{Prop_G}, 
for a constant $C > 0$ to be determined later on, we have 
\begin{equation}
\begin{aligned}
& \mathbb{E}\Bigg(\sqrt{\dfrac{\alpha^2 \|{\bf G}^\dagger\|_2^2}{1 + \beta^2 
\|{\bf G}^\dagger\|_2^2}}\Bigg) = \mathbb{E}(g(\|{\bf G}^\dagger\|_2)) \\ & =\int_{0}^{\infty} g\prime(x)\mathbb{P}\{\|{\bf G}^\dagger\|_2 \ge x\}dx \\
& \le \int_{0}^{C} g\prime(x)dx +
\int_{C}^{\infty} g\prime(x)\mathbb{P}\{\|{\bf G}^\dagger\|_2 \ge x\}dx \\
& \le \sqrt{\dfrac{\alpha^2 C^2}{1 + \beta^2 C^2}} + 
\int_{C}^{\infty} \dfrac{\alpha^2
\Big({\dfrac{p+1}{\text{e}\sqrt{\ell}}x\Big)^{-(p+1)}}}{(1 + \beta^2x^2)^2\sqrt{\dfrac{\alpha^2 x^2}{1 + \beta^2x^2}}}dx\\
& = \sqrt{\dfrac{\alpha^2 C^2}{1 + \beta^2 C^2}} + 
\dfrac{\alpha^2C^2}{(p-1)(1 + \beta^2 C^2)^2\sqrt{\dfrac{\alpha^2 C^2}{1 + \beta^2 C^2}}} \times \\ 
& \quad \underbrace{\int_{C}^{\infty}\Big({\dfrac{p+1}{\text{e}\sqrt{\ell}}x\Big)^{-(p+1)}}}_{\Big({\dfrac{p+1}{\text{e}\sqrt{\ell}}C\Big)^{-(p+1)}}}.
\end{aligned}
\notag
\end{equation}

Likewise, we seek $\nu>0$ such that

\begin{equation}
\begin{aligned}
& \sqrt{\dfrac{\alpha^2 C^2}{1 + \beta^2 C^2}} + 
\dfrac{\alpha^2C^2\Big({\dfrac{p+1}{\text{e}\sqrt{\ell}}C\Big)^{-(p+1)}}}{(p-1)(1 + \beta^2 C^2)^2\sqrt{\dfrac{\alpha^2 C^2}{1 + \beta^2 C^2}}} \\ 
& \le \sqrt{\dfrac{\alpha^2 \nu^2}{1 + \beta^2 \nu^2}}.
\end{aligned}
\notag
\end{equation}

The solution satisfies
\begin{equation}
\begin{aligned}
\nu\ge \dfrac{C^2}{\sqrt{C^2 - \dfrac{C^2}{p-1}\Big({\dfrac{p+1}{\text{e}\sqrt{\ell}}C\Big)^{-(p+1)}}}}.
\end{aligned}
\notag
\end{equation}

The value $\nu = \dfrac{4\text{e}\sqrt{\ell}}{p+1}$ satisfies this inequality for $C = \Big({\dfrac{\text{e}\sqrt{\ell}}{p+1}\Big)\Big({\dfrac{2p}{p-1}\Big)^{1/(p+1)}}}$. \QEDB


%


\bibliographystyle{IEEEtran}
\bibliography{mybibfile}
%








\end{document}